\documentclass[12pt]{article}
\usepackage[english]{babel}
\usepackage{amsmath}
\usepackage{amssymb}
\usepackage{amsfonts}
\usepackage{amsthm}
\usepackage{graphics,graphicx}
\usepackage{color}

\usepackage{graphicx,amsmath,amsthm,amsrefs}

\usepackage{color}
\usepackage{appendix}
\definecolor{black}{gray}{0}
\definecolor{ttgray}{gray}{0.5}
\definecolor{bfred}{rgb}{0.4,0,0}
\definecolor{emgreen}{rgb}{0,0.3,0}
\definecolor{emmagenta}{rgb}{0.6,0.07,0.07}
\definecolor{sfgb}{rgb}{0,0.3,0.3}
\definecolor{mathblue}{rgb}{0,0,0.4}

\newcommand{\chegakh}[1]{\{\kern -.25em [#1]\kern -.25em \}}
\newcommand{\Dchegakh}[1]{\{\kern -.3em\{\kern -.25em [#1]\kern -.25em\}\kern -.3em\}}

\newtheorem{Question}{Question}

        \newtheorem{thm}{Theorem}[section]

    \PassOptionsToPackage{british}{babel}
    
    \newtheorem{lem}[thm]{Lemma}

    \theoremstyle{definition}

    \theoremstyle{remark}
    \newtheorem{rem}[thm]{Remark}

\begin{document}

 \begin{center}

\vspace{5mm}

{\large\bf {The Wiener--Hopf Technique, its Generalisations and
    Applications: Constructive and Approximate Methods}}

\vspace{10mm}
{\bf  Anastasia KISIL$^1$,David ABRAHAMS$^2$ Gennady MISHURIS$^3$, Sergei ROGOSIN$^{4}$}

\vspace{3mm}
{\footnotesize\it $^1$ Department of Mathematics, The University of Manchester,
Manchester, M13 9PL, UK;

\vspace{1mm}
e-mails: anastasia.kisil@manchester.ac.uk}

\vspace{3mm}
{\footnotesize\it $^2$Isaac Newton Institute for Mathematical Sciences
University of Cambridge, 20 Clarkson Road, CB3 0EH Cambridge, UK;

\vspace{1mm}
e-mails: i.d.abrahams@newton.ac.uk}

\vspace{3mm}
{\footnotesize\it $^3$Aberystwyth University, Penglais, SY23 3BZ Aberystwyth, UK;

\vspace{1mm}
e-mails: ggm@aber.ac.uk}

\vspace{3mm}
{\footnotesize\it $^4$Belarusian State University, Nezavisimosti Ave., 4, 220030 Minsk, Belarus;

\vspace{1mm}
e-mail: rogosin@bsu.by}

\end{center}

{\footnotesize {\bf Abstract.} This paper reviews
  the modern state of the Wiener--Hopf factorization method and its
  generalizations. The main constructive results for matrix Wiener--Hopf are presented,
  approximation methods are outlined and the main areas of
  applications are mentioned. The aim of the paper is to sketched some perspective of the development of this method, importance of bringing together pure and applied analysis to most effectively use of the Wiener--Hopf technique.

  \vspace{1mm} {\bf Key words and phrases:} Wiener--Hopf,
Riemann--Hilbert,  factorization, partial indexes, Riemann boundary value problem, applications

\vspace{1mm} {\bf AMS Subject Classification 2010:} primary 30E25, 45E10, 47A68, 47B35; secondary 74J20,
74N10, 78A45}

\begin{abstract}
This paper reviews the modern state of the Wiener--Hopf factorisation method and
its generalizations. The main constructive results for matrix Wiener--Hopf problems
are presented, approximate methods are outlined and the main areas of applications
are mentioned. The aim of the paper is to offer an overview of the development of
this method, and demonstrate the importance of bringing together pure and applied
analysis to most effectively employ the Wiener--Hopf technique.
\end{abstract}


\section{Introduction}
\label{sec:into}

The Wiener--Hopf method has been motivated by interdisciplinary interests ever since its
inception. It resulted from a collaboration between Norbert Wiener, who
worked on stochastic processes, and Eberhard Hopf, who worked on
partial differential equations (PDEs). The method was first described
in their joint article \cite{WieHop31} where they study a
convolution type integral equation for \(f(x)\) on a semi-axis
\begin{equation}
\label{eq:WH1}
\int\limits_{0}^{\infty} k(x - t) f(t) dt = g(x) \quad x > 0,
\end{equation}
 with \(k(x)\) and \(g(x)\) given.

The study of this equation was principally motivated by an interest of one of the
authors, in a differential equation governing the radiation equilibrium
of stars~\cite{Hop34}. The Wiener--Hopf equation arises from extending the integral equation \eqref{eq:WH1}
into \(x<0,\)
\begin{equation}
\label{eq:WH1_ex}
\int\limits_{0}^{\infty} k(x - t) f(t) dt =
\begin{cases}
   g(x)  & \text{for }  x > 0,\\
   h(x)     & \text{for } x \le 0,
  \end{cases}
\end{equation}
for some function \(h(x)\). Note that, although \( h(x) \) is unknown, it is not independent as it uniquely defined by the
left hand side of (\ref{eq:WH1_ex}) once \( f(x) \) is determined. Applying the Fourier transform to
\eqref{eq:WH1_ex} results in an equation of type \eqref{eq:WH2} (the
interaction of a Fourier transform with its convolution has to be used, and
analytic properties of half-range Fourier transforms employed\footnote{In fact, in the
original paper \cite{WieHop31} the Laplace
transform  was applied and the corresponding functional equation is
valid in a vertical strip. Interestingly the original paper
contained a mistake; it  claimed that \(\gamma(z)^n=(z^2-1)^{n/2}\to z^n\) as
\(z \to \pm \infty\)
where \(n \) was the number of zeros of the kernel in the strip of
analyticity. This was pointed out in~\cite[Chapter~16]{Ga-Che}), and
resulted from an incorrect choice of branch cut. However, it did
not influence the actual results they obtained in the second
part of the paper, since they assumed that the kernel is even and hence
\(n\) is even.}). This laid the foundation for the study of scalar
Wiener--Hopf equations.

Independently, the theory of the more general Riemann--Hilbert boundary value
problem has been developed. The Riemann--Hilbert boundary value problem for a particular case was first
formulated by Riemann as a part of the problem of construction of
complex differential equation of the Fuchs type with a given
monodromy. In the latter form, the problem was presented by Hilbert
as his 21st mathematical problem for the 20th century (see
\cite{Bol90}). Plemelj \cite{Ple08} proposed a method of solution to
the Riemann--Hilbert boundary value problem based on the reduction to a case of the homogeneous Riemann boundary value problem proving its solvability by means of the Fredholm alternative\footnote{For a long time it was supposed that Plemelj's work gave a complete answer on the Hilbert question. But in the late 1980s Bolibrukh (see, e.g., \cite{Bol90}) showed that
the proof of Plemelj is incomplete and that the negative answer is also possible. In fact, it was shown that Plemelj's positive
answer for the question is valid for the so-called ``regular'' (see,~\cite[p. 7]{Bol90}) variant of the Riemann--Hilbert boundary value problem.}. An intensive study of the Riemann--Hilbert boundary value problem began in the 1930s. Properties of the factorisation problem and its role in the study of the Wiener--Hopf integral equation were outlined in the seminal paper \cite{GohKre58}. Further development of the factorisation theory was described in the book by Litvinchuk and Spitkovskii \cite{Spit_book}.
Several known facts on factorisation were presented in a recent survey
\cite{RogMis16} (see also \cite{Spit}), where one can find, in
particular, an extended list of references on the subject.

Simultaneously, the general theory of integral equations was developed
intensively at the beginning of 20th century. Of note are the
contributions by Hilbert \cite{Hil12}, Fredholm \cite{Fre03}, Volterra
\cite{Vol09} and Carleman \cite{Car22}. The latter two authors also considered
integral equations with kernels depending on the difference of arguments,
but containing integrals with variable limits or along a finite
interval, respectively. Nevertheless, the Wiener--Hopf
integral equation \eqref{eq:WH1} became an important independent
subject of research due to its wide range of applicability and many
connections. The Wiener--Hopf method has enabled researchers to analytically
solve numerous previously intractable integral and partial differential
equations, as well as many boundary value problems.
To date, it still remains the standard method for solving a wide class of
canonical physical problems. The
method is an important cornerstone (amongst other methods) in areas of
pure mathematics, for example in development of the abstract theory of
singular linear operators, pseudo-differential operators etc.
Additionally, the exact solutions  form the basis
of approximate techniques applicable to more complicated problems. Moreover, 
even in case when the constructive (analytical) construction is available, it is not always possible to 
realise numerically with a confidence \cite{AAM_2021}.

Recently, the Isaac Newton Institute for Mathematical Sciences, Cambridge UK, organised a
one-month research programme ``Bringing pure and applied analysis
together via the Wiener--Hopf technique, its generalisations and
applications'', where the most important results in the area were
presented and further developments of the method were considered~\cite{NSR_20}. This
review was inspired by that programme. The aim of this paper is to show
the beauty, principal features and perspectives of the Wiener--Hopf technique.
The theory of scalar Wiener--Hopf equations is now very rich and well
developed. In contrast, much less is known about matrix Wiener--Hopf
equations. These are a natural extension of the scalar case, and
enable us to model more advanced problems derived from applications~\cite{history}.
As it stands, solutions to matrix Wiener--Hopf problems have to be constructed on a case-by-case
basis, or in an approximate fashion; the main approaches will be
outlined in this review.
The Wiener--Hopf technique is currently used in a wide range of
disciplines including acoustics, finance, L\'{e}vi processes,
hydrodynamics, elasticity, potential theory and electromagnetism.
As such, there are several disjoint
communities who rely on the Wiener--Hopf method, and developments in
one field sometimes go unnoticed in another. It is intended that this
review helps to remove this artificial intradisciplinary and
interdisciplinary boundary.

In summary, the Wiener--Hopf technique is captivating for many reasons:
\begin{itemize}
\item It enables us to solve numerous physical problems motivated by real world
applications modelled, for example, by partial differential equations and stochastic processes.
\item Subtle and diverse analytical methods are required to obtain
solutions.
\item Solutions of the Wiener--Hopf problems are remarkably revealing, thus often allowing for explicit evaluation of its asymptotic forms at all singular points and at infinity.
\item Computing solutions to obtain practically useful results often requires involved and innovative numerical analysis.
\end{itemize}

This review offers a short introduction to the application, theory and numerical implementation of the Wiener--Hopf method; each element
forms an important part of the whole.

The review has the following structure. The following section summarises results
related to the Wiener--Hopf technique, related equations and applications. Sections~\ref{constructive} lists the main
constructive methods for certain classes of matrix Wiener--Hopf equations allowing for explicit factorisation.
Section~\ref{sec:approx} outlines several approximate methods available for the cases where exact constructive
methods do not apply. The last section draws conclusions and discusses further directions.

\section{Preliminaries}
\label{sec:pre}

\subsection{The Wiener--Hopf problem}
\label{sec:W-H}

We begin by formulating the basic Wiener--Hopf equation and discussing its
standard method of its solution~\cite{bookWH}. This will be compared and contrasted with the
Riemann--Hilbert boundary value problem at the end of the next
section. The Wiener--Hopf equation has a simpler form and
historically was the first to be discovered, so forms a natural
starting point.

Define {\it upper half-plane} ${\mathcal H}^{+} =
\{{\mathrm{Im}}\, z > a; a < 0\}$, {\it lower half-plane} ${\mathcal H}^{-} =
\{{\mathrm{Im}}\, z < b; b > 0\}$ and {\it the Wiener--Hopf strip} ${\mathcal H} =
{\mathcal H}^{+}\cap {\mathcal H}^{-} = \{a < {\mathrm{Im}}\, z < b\}$
($a$ and $b$ are real numbers, typically small). Then the classic
{\it Wiener--Hopf problem} is, for given $K(z)$ and $C(z)$ analytic in
${\mathcal H}$, and
\begin{equation}
\label{eq:WH2}
K(z)\Phi_+(z)+ \Psi_-(z)+C(z)=0, \quad z\in {\mathcal H},
\end{equation}
to find $\Phi_{+}(z)$ an unknown function analytic in ${\mathcal H}^{+}$  and  $\Psi_{-}(z)$  an unknown function
analytic in ${\mathcal H}^{-}$. The function $K(z)$ is called {\it the Wiener--Hopf kernel corresponding to
  \eqref{eq:WH2}}.

For the sake of completeness we are going to outline the usual
procedure used to solve the Wiener--Hopf problem in the case of
canonical scalar factorisation.
The key  step in the Wiener--Hopf procedure is to find \(K_+\)
and \(K_-\) (\emph{the Wiener--Hopf factors}) such that:
\begin{equation}
\label{WH_fact}
K(z)=K_+(z)K_-(z),\quad  z\in {\mathcal H}.
\end{equation}
with \(K_\pm\) and their inverses analytic in ${\mathcal H}^{\pm}$.
This product is referred to as a \emph{factorisation of the
kernel} and is the most important part of the procedure; it is straightforward to accomplish for scalar kernels.

Multiplying through by \(K_-^{-1}\) we obtain:
\[K_+ \Phi_+ ( z)+ K_-^{-1}  \Psi_-(z)+ K_-^{-1}C(z)=0,\quad z\in {\mathcal H}.\]
Now performing an \emph{additive splitting} we have:
\[ K_-^{-1}C(z)=C_-(z)+C_+(z),\quad z\in {\mathcal H},\]
where \( C_\pm(z) \) are analytic in their indicated half-planes \({\mathcal H}^\pm\), respectively.
Finally rearranging the equation yields
\begin{equation}
\label{WH_fact_estimates}
K_+ \Phi_+ ( z) +C_+( z)=-K_-^{-1}
\Psi_-(z)-C_-(z), \quad  z\in {\mathcal H}.
\end{equation}
The right hand side of \eqref{WH_fact_estimates} is a function
analytic in ${\mathcal H}^{+}$ and the left hand side of of
\eqref{WH_fact_estimates} is analytic in ${\mathcal H}^{-}$. Hence, together they offer
analytic continuation from ${\mathcal H}$ into the whole of the $z$-plane, and so
each side is equal to an entire function $J(z)$, say. We also require \emph{at most
algebraic growth} of each side of \ref{WH_fact_estimates} as $|z| \to \infty$ in the
respective half-planes of analyticity, i.e.\
\[J(z)\leq {\cal O} (|z|^{n}), \quad |z| \to \infty,\]
for some constant $n$. Applying the extended form of Liouville's theorem implies that
\(J(z)\) is a polynomial of degree less than or equal to
\(n\). Hence \(\Phi_+ ( z)\) and \(\Psi_-(z)\) are
determined up to \(\lfloor n \rfloor +1\) unknown constants (typically \(n\leq 0\) or
\(1\), and the constants can be fixed using extra information about the behaviour of the solution).

If  $\Phi_{+}(z)$, \(\Psi_-(z)\) and \(C(z)\) in \eqref{eq:WH2} are vectors of
functions with a matrix of functions \(K(z)\), then~\eqref{eq:WH2} is
called  \emph{matrix/vectorial Wiener--Hopf equations}. All the above steps in
the solution are analogous. The difficulty lies in performing the
 multiplicative factorisation (\ref{WH_fact}), which becomes the main challenge to
 overcome and is addressed in Sections~\ref{constructive} and~\ref{sec:approx}.

\subsection{Riemann--Hilbert boundary value problem}
\label{sec:R-H}

The Wiener--Hopf equation is closely related to the Riemann--Hilbert boundary value
problem\footnote{It would be more historically accurate to call it the
  Riemann boundary value problem, as still used by Gakhov's school; however, we will use the title Riemann--Hilbert boundary
  value problem, introduced by Muskhelishvili, as it is more common nowadays.}~\cite{Gak77}. The latter problem is to
determine two functions or vectors of functions  $\Phi^{\pm}(z)$,
which are analytic in the complementary domains ${\mathcal D}^{\pm}$
(i.e. ${\mathcal D}^{+}\cup {\mathcal L} \cup {\mathcal D}^{-} =
{\mathbb C}$, with a given simple closed curve ${\mathcal L} = \partial {\mathcal D}^{+} = \partial{\mathcal D}^{-}$) and satisfy on ${\mathcal L}$ the following linear condition
\begin{equation}
\label{eq:R1}
\Phi^{+}(t) = G(t) \Phi^{-}(t) + g(t), \quad t\in {\mathcal L},
\end{equation}
with the function/matrix $G(t)$ and the function/vector $g(t)$ given on ${\mathcal L}$.

In the scalar case the complete solution of the Riemann problem \eqref{eq:R1} was given by Gakhov (see \cite{Gak77}). He has shown that solvability of the problem depends on the {\it Cauchy index} (winding number) of the coefficient $G$
\begin{equation}
\label{R2}
\kappa = {\mathrm{ind}}_{\mathcal L} G = \frac{1}{2\pi} \int\limits_{\mathcal L} d (\arg G(t)).
\end{equation}
In the case of H\"older-continuous functions $G$ and $g$, a simple
smooth contour ${\mathcal L}$ and in the case of non-negative index
$\kappa\geq 0$ the solution of \eqref{eq:R1} has the
form\footnote{These formulas are valid for ${\mathcal L} = {\mathbb
    R}$ and under weaker assumptions on $G, g$ (see
  \cite{Gak77}). }
\begin{equation}
\label{eq:R3}
\Phi^{\pm}(z) = X^{\pm}(z) \left[\Psi^{\pm}(z) + P_{\kappa}(z)\right],\; z\in {\mathcal D}^{\pm},
\end{equation}
where $P_{\kappa}(z)$ is a polynomial of order $\kappa$, and
\begin{equation}
\label{eq:R4}
X^{+}(z) = \exp \{\Gamma^{+}(z)\},\; X^{-}(z) = \exp \{z^{-\kappa} \Gamma^{-}(z)\},
\end{equation}
\begin{equation}
\label{eq:R5}
\Gamma^{\pm}(z) = \frac{1}{2\pi i} \int\limits_{\mathcal L} \frac{\log \left[t^{-\kappa} G(t)\right] dt}{t - z}, \; z\in {\mathcal D}^{\pm},
\end{equation}
\begin{equation}
\label{eq:R6}
\Psi^{\pm}(z) = \frac{1}{2\pi i} \int\limits_{\mathcal L} \frac{g(t)}{X^{+}(t)}\frac{dt}{t - z}, \; z\in {\mathcal D}^{\pm}.
\end{equation}
In the case of negative index $\kappa$ the unique solution has the
form $\Phi^{\pm}(z) = X^{\pm}(z) \Psi^{\pm}(z)$ provided that a finite
number of solvability conditions
$$
\int\limits_{\mathcal L} \frac{g(t)}{X^{+}(t)} t^{k - 1} dt = 0, \; k = 1, 2, \ldots, - \kappa - 1,
$$
are satisfied (otherwise there is no solution).

The scalar and matrix function factorisation, related to the Riemann--Hilbert boundary value problem, was introduced by Birkhoff \cite{Bir03} and also in a more straightforward
form by Gakhov \cite{Gak50}. Factorisation means obtaining the following
representation of an $n\times n$ matrix function $G$, given on a bounded
curve ${\mathcal L}$ (curve such that $0\in {\mathcal D}^{+}$, $\infty\in {\mathcal D}^{-}$), as
\begin{equation}
\label{eq:F1}
G(t) = G^{+}(t) \Lambda(t) G^{-}(t), \; t\in {\mathcal L},
\end{equation}
where $G^{\pm}(t)$ are boundary values of bounded analytic zero-free matrices $G^{\pm}(z), z\in {\mathcal D}^{\pm}$, and
\begin{equation}
\label{eq:F2}
\Lambda(t) = {\mathrm{diag}} \{t^{\kappa_1},t^{\kappa_2},\ldots,t^{\kappa_n}\},
\end{equation}
where $\kappa_j\in {\mathbb Z}, j=1,2,\ldots,n,$ are integer numbers
called {\it partial indices}. Representation \eqref{eq:F1} is called a {\it
  left-sided factorisation}. Interchanging $G^{+}$ and $G^{-}$ we
arrive at the {\it right-sided factorisation}. Respectively, the partial indices are called left (right) partial indices.
In general, they are not the same, however their sum is always
equal to the index of \(\det G\). In the case ${\mathcal  L} = {\mathbb R}$ the variable $t$ in \eqref{eq:F2} is replaced by the
ratio $\left(\frac{t-i}{t+i}\right)$ (similar changes are need in
\eqref{eq:R4} and \eqref{eq:R5}).
In this case one needs to assume the continuity of $G$ at infinity.

Unfortunately, in the general case, there exists no method of determining the partial
indices, while such information is important to reconstruct the factors, $G^{\pm}$, and is
instrumental in proving the stability of any respective algorithm.
Those issues are discussed further in Section~\ref{sec:approx}. Factorisation, i.e. determination of the partial indices  $\kappa_j\in {\mathbb Z},
j=1,2,\ldots,n,$ and factors $G^{+}, G^{-}$ is an independent problem
playing a crucial role in many theoretical and applied problems.

We would like to make here a few comments on formal differences between the Wiener--Hopf and
Riemann--Hilbert problem. Note that \eqref{eq:F1} is different to
\eqref{WH_fact} since equation \eqref{WH_fact} does not contain \(\Lambda(t)\). It is
because it is common in the Wiener--Hopf literature~\cite{bookWH} to implicitly assume that
\(\kappa=0\) in formulae \eqref{eq:R4} and
\eqref{eq:R5}, which is often the case in applications. On the other hand, \eqref{eq:F1} can be always
presented in the form \eqref{WH_fact} by further factorising the matrix $\Lambda(t)=\Lambda^+(t)\Lambda^-(t)$, and then redistributing the respective diagonal terms.
However, this destroys some of the assumptions on the behaviour of the factors $G^{\pm}$: the behaviour at infinity, satisfaction of the invertibility property (i.e.\ it introduces zeros into the half-planes of analyticity) etc.
Thus, this form of the factorisation is tightly linked to the
conditions imposed on the factors (which space of analytic function is
chosen)~\cite{Spit_book}. Note also that the factors \(K_\pm\)
themselves, \eqref{WH_fact}, would in general have non-zero partial indices.
Another important difference between the approaches is that \eqref{WH_fact} is valid on a strip (and thus, knowledge of the possible poles and zeros inside the strip provides
useful information for solving the problem under consideration), while the factorisation
\eqref{eq:F1} is defined on a line (that can be analytically extended if
possible). These differences are also discussed in~\cite{MyStrip}. Finally, we would like to stress that there are authors associating the Wiener-Hopf problem with the Riemann-Hilbert one and vice versa that only underlines the fact that the division is rather artificial and depends on the research communities. To avoid any unintentional misleading, we discuss particular factorisation techniques below as they have appeared in the original presentations, with minor comments where necessary.

\subsection{Canonical matrix}
\label{subsec:can_mat}

The factorisation problem \eqref{eq:F1} is almost equivalent to the
homogeneous
(i.e. with $g(t)\equiv 0$) Riemann--Hilbert boundary value problem \eqref{eq:R1}:
\begin{equation}
\label{R0}
\Phi^{+}(t) = G(t) \Phi^{-}(t), \; t\in {\mathcal L}.
\end{equation}
In \cite{Mus68} the conditions for the existence of a solution to (\ref{R0}) are described.
If the solution to (\ref{R0}) exists one can find\footnote{This is a
  crucial problem, highly dependent on the properties of the matrix
  coefficient $G$, and NOT always possible to be performed in a constructive way.} a so called {\it fundamental system of solutions}  $\Phi_1(z), \ldots, \Phi_n(z)$ and the corresponding {\it fundamental matrix} $\Phi(z) = \left(\Phi_1(z), \ldots,
\Phi_n(z)\right)$ with vectors $\Phi_j(z)$ being their columns. Fundamentality of the system/matrix means that the
determinant of the fundamental matrix cannot be identically zero.

The next step (which is more constructive, see \cite[p. 520]{Mus68}) is to
transform the fundamental system of solutions to a normal form. By
definition  the {\it normal system of solutions}
$\Psi_1(z),\ldots,\Psi_n(z)$ to  (\ref{R0}) is the fundamental system
such that the determinant of the fundamental matrix does not vanish anywhere on ${\mathbb C}$ (including the curve ${\mathcal L}$). The matrix $\Psi(z) = \left(\Psi_1(z) \ldots \Psi_n(z)\right)$ of the normal system is called the {\it normal matrix}.

By elementary transformations the normal system of solutions can be
transformed to the so called {\it canonical system}. The canonical system
is usually denoted $X_1(z), \ldots, X_n(z)$ and its matrix is $X(z) =
\left(X_1(z) \ldots X_n(z)\right)$.  The properties of the canonical
system are
\begin{enumerate}
\item The canonical system is the normal system of solutions to
  (\ref{R0}), i.e. the determinant of the canonical matrix $\det\,
  X(z)$ is nowhere vanishing in ${\mathbb C}$.\\
\item  The sum of the orders $\kappa_j$  at infinity of the columns of the matrix $X  (z)$ (i.e. solutions $X_j(z)$) is equal to the Cauchy index of the determinant $\det X(z)$.
\end{enumerate}
They are called partial indices of the Riemann--Hilbert boundary value problem
(\ref{R0}).
The integer numbers $\kappa_j$ are
the same as partial indices of the
factorisation \eqref{eq:F1} of the matrix function $G(t)$ (the
coefficient in  (\ref{R0})).

Whenever the canonical matrix is constructed, the solution of factorisation problem \eqref{eq:F1} is presented in the form
\begin{equation}
\label{bvp1}
G^{+} = X^{+}, \;\;\; G^{-} = \Lambda^{-1} (X^{-})^{-1},
\end{equation}
where $\Lambda(z) = \textrm{diag} \left\{z^{\kappa_1}, \ldots,
  z^{\kappa_n}\right\}$ (this applies when \(0 \in{\mathcal D}^{+}\)).

Note that for a canonical matrix it is said
that the matrix $X^{-}(z)$ has {\it the normal form at infinity} for a
factorisation \eqref{eq:F1}.
An example of the use of this method is presented in
Section~\ref{constructive}\ref{triangle}.

\subsection{Other convolution type integral equations}

Related to the integral equations of the convolution type~\eqref{eq:WH1} on the half-axis are paired equations
\begin{equation}
\label{dual}
\begin{matrix}
\displaystyle
f(x)+\int\limits_{-\infty}^{\infty} k_1(x - t) f(t) dt&=g(x),
\quad x>0,\\[2mm]
\displaystyle
f(x)+\int\limits_{-\infty}^{\infty} k_2(x - t) f(t) dt&=g(x), \quad x<0,
\end{matrix}
\end{equation}
and the so-called transpose to the dual equations:
\begin{equation}
\label{dual_transpose}
f(x)+\int\limits_{0}^{\infty} k_1(x - t) f(t) dt+\int\limits_{-\infty}^{0} k_2(x - t) f(t) dt=g(x), \quad x\in \mathbb{R}.
\end{equation}
They can be both solved in a very similar manner to~\eqref{eq:WH1} by
extending the range, taking the Fourier transform of each of the
equations, and eliminating one unknown, to give rise to the
usual Wiener--Hopf equation
\begin{equation}
\label{dual_WH}
F^+(s)= R(s)F^-(s)+T(s),\quad R(s)=\frac{1+K_2(s)}{1+K_1(s)},
\end{equation}
where $K_j(s)$ are Fourier transforms of the kernels $k_j(x)$, $j=1,2$, and $F^\pm$, $T(s)$ are different for each equation
(details can be found in~\cite{FGK94, Gak77}).

\subsection{Carleman type boundary value problem}

A closely related equation, but much less known is
the so called smooth-transition equation
\begin{equation}
\label{eq:Car}
f(x)+\int\limits_{-\infty}^{\infty} k_1(x - t) f(t) dt- g(x)+
e^{-x}\left\{ f(x)+\int\limits_{-\infty}^{\infty} k_2(x - t) f(t) dt- g(x)\right\}=0,
\end{equation}
for \(-\infty \le x \le \infty \). The interesting feature of this
equation is that for \(x\) large and positive it is close to
\begin{equation}
\label{eq:Car_small}
f(x)+\int\limits_{-\infty}^{\infty} k_1(x - t) f(t) dt- g(x)=0, \quad
x \gg1,
\end{equation}
 And for \(x\) large and negative it is close to
\begin{equation}
\label{eq:Car_large}
f(x)+\int\limits_{-\infty}^{\infty} k_2(x - t) f(t) dt- g(x)=0, \quad x\ll-1.
\end{equation}
 which is similar to \eqref{dual}.
Finally, for values of $x=0$ close to zero, the equation takes the form:
\begin{equation}
\label{eq:Car_middle}
2f(x)+\int\limits_{-\infty}^{\infty} k_1(x - t) f(t) dt+
\int\limits_{-\infty}^{\infty} k_2(x - t) f(t) dt- g(x)=0, \quad |x|\ll1.
\end{equation}
This analysis explains the name of the equation; it behaves as the dual equation when $|x|\gg1$ and as the transpose to the dual when $|x|\ll1$.

In order to analyse this equation, define the following~\cite{Ga-Che}
\begin{equation}
\label{eq:Car1}
\phi(x)=f(x)+\int\limits_{-\infty}^{\infty} k_2(x - t) f(t) dt- g(x), \quad -\infty \le x \le \infty.
\end{equation}
Then the original equation \eqref{eq:Car} can then be written as
\begin{equation}
\label{eq:Car2}
f(x)+\int\limits_{-\infty}^{\infty} k_1(x - t) f(t) dt- g(x)+
e^{-x}\phi(x)=0.
\end{equation}
After application of Fourier transform to \eqref{eq:Car1}
and \eqref{eq:Car2} (again, with a capital letter indicating
the transform of a lower-case letter) we arrive at the following
equations
\begin{equation}
  \begin{aligned}
  \Phi(z)&=&(1+K_2(z))F(z)-G(z),\quad -\infty \le z \le \infty,\\
  -\Phi(z+i)&=&(1+K_1(z))F(z)-G(z), \quad -\infty \le z \le \infty.
\end{aligned}
  \label{eq:pp}
\end{equation}
Eliminating \(F(z)\) we arrive at
\begin{equation}
  (1+K_1(z))\Phi(z)=-(1+K_2(z))\Phi(z+i)+ (K_2(z)-K_1(z))G(z), \quad -\infty \le z \le \infty.
\end{equation}
This is not a standard form for a Riemann--Hilbert boundary value problem but it
can be transformed to one by defining a new function
\begin{equation}
  \label{eq:trans}
  w(\psi)=\frac{1}{\sqrt{\psi}}\,\Phi\left(\frac{\ln(\psi)}{2\pi}\right),
\end{equation}
then the jump in the Riemann--Hilbert boundary value problem is across the cut
\(\arg{\psi}=0\).

The smooth transition equation \eqref{eq:Car} is a special case of a more
general Carleman type boundary value problem: to find \(\Phi(z)\)
analytical in a domain \({\mathcal D}\)
and satisfy on part of the boundary of \({\mathcal D}\) denoted ${\mathcal S}$ the following condition
\begin{equation}
\label{eq:R1c}
\Phi(t) = E(t) \Phi(m(t)) + g(t), \quad t\in {\mathcal S},
\end{equation}
with the functions $E(t)$, \(m(t)\) and  $g(t)$
given on ${\mathcal S}$~\cite{Ga-Che}*{\S~15.1}. Also it is required that \(m(t)\) is
continuous and maps to the rest of the boundary of \({\mathcal
  D}\). Additionally we require that \(t\) and \(m(t)\) transverse the boundary in different
directions. Sometimes its solution can be reduced to a
Riemann--Hilbert boundary value problem, but in general no closed form solution is known~\cite{Lit_shift}. Many applied problems are reduced to certain types of
 functional-difference equation (related to Wiener-Hopf and
 Riemann-Hilbert problems). We mention here the problems of
 applied mechanics~\cite{Atkinson_94,Antipov03, Mishris96,
   Mishuris99,Mishuris01} and
 diffraction theory~\cite{AntSil04, Budaev95, Croisille99}.

\subsection{Discrete Wiener--Hopf equation}

In this section we will look at an infinite system of equations
of Wiener--Hopf type, a discrete analogue of integral equation
  \eqref{eq:WH1}. For all
infinite sequence \(a_n\) considered here we will assume
the following holds for some constant \(M\)
\begin{equation}
  \label{eq:DHol}
  |a_n|<\frac{M}{n^{1+\lambda}}, \quad 0<\lambda<1.
\end{equation}
We will also assume that the Fourier series
\(A(\theta)=\sum_{n=-\infty}^{\infty} a_n e^{in\theta}\) is H\"{o}lder
continuous.

The \emph{discrete  Wiener--Hopf} equations for \(x_n\) is
\begin{equation}
  \label{eq:DW-H}
  \sum_{k=0}^{\infty} a_{n-k}x_k=c_n, \quad n=0,1,2 \dots,
\end{equation}
where the infinite sequences \(a_n\) and \(c_n\) are given.
This is called an infinite system with a Toeplitz operator due to a
presence of \(n-k\) in the index of \(a\). Similar to the integral
equation case, we first need to extend this equation for \(n\) negative
\begin{equation}
  \label{eq:DW-H1}
  \sum_{k=0}^{\infty} a_{n-k}x_k=\begin{cases}
   c_n  & \text{for } n=0,1,2 \dots \\
   d_n      & \text{for } n=-1,-2 \dots
  \end{cases}
\end{equation}
for some \(d_n\), which are unknown for \(n=-1,-2 \dots\). Next we apply the discrete
analogue of Fourier transform, the \(Z\)-transform. The \(Z\)-transform of
a sequence \(a_n\) is defined as
\begin{equation}
  \label{eq:ztran}
  A(t)=\sum_{k=-\infty}^{\infty} t^ka_k,
\end{equation}
and as before is denoted by a capital letter. We multiply the \(n\)th
equation in \eqref{eq:DW-H1} by \(t^n\) and sum over all equations:
\begin{equation}
  \label{eq:Dsum}
   \sum_{n=-\infty}^{\infty}t^n \sum_{k=0}^{\infty} a_{n-k}x_k= \sum_{n=0}^{\infty}t^nc_n+\sum_{n=-1}^{-\infty}t^nd_n.
 \end{equation}
The main property is the interaction of the \(Z\)-transform with a Toplitz
operator
\begin{equation}
  \label{eq:Dsum_pro}
   \sum_{n=-\infty}^{\infty}t^n \sum_{k=0}^{\infty} a_{n-k}x_k=\sum_{k=0}^{\infty} t^kx_k \sum_{n=-\infty}^{\infty}t^{n-k} a_{n-k}
   =\sum_{k=0}^{\infty} t^k x_k \sum_{m=-\infty}^{\infty}t^{m} a_{m},
 \end{equation}
which corresponds to the Fourier transform changing a convolution
into a product. Thus, via a \(Z\)-transform, equation \eqref{eq:DW-H1} becomes
\begin{equation}
  \label{eq:CW-H}
  A(t)X^+(t)=C(t)+D^-(t),
\end{equation}
which holds on the unit circle, and the superscript \(+/-\)
indicates the function is analytic inside/outside the unit
circle. This Riemann--Hilbert boundary value equation can be solved as before to find \(X^+(t)\) and then the inverse
of the \(Z\)-transform applied to recover \(x_n\).

A slightly more general infinite system which can be solved in an identical fashion is the discrete version of equation \eqref{dual}
\begin{equation}
  \label{eq:DW-dual}
\begin{matrix}
\displaystyle
  \sum_{k=-\infty}^{\infty} a_{n-k}x_k=
   c_n,  & \text{for } n=0,1,2 \dots \\[2mm]
\displaystyle
  \sum_{k=-\infty}^{\infty} b_{n-k}x_k=
   d_n,      & \text{for } n=-1,-2 \dots
\end{matrix}
\end{equation}
where $a_k$, $b_k$, $c_k$ and $d_k$ are known.
Another version of the discrete equation is the analogue to the transpose to the dual equation \eqref{dual_transpose}:
\begin{equation}
  \sum_{k=0}^{\infty} a_{n-k}x_k+\sum_{k=-1}^{-\infty} b_{n-k}x_k=c_n,
  \quad n=0, \pm 1,\pm 2 \dots
\end{equation}

Another interesting discrete system which has an analytic solution via a
Wiener--Hopf technique is
\begin{equation}
  \sum_{k=-\infty}^{\infty} a_{n-k}x_k-c_n+r^n\left(
    \sum_{k=-\infty}^{\infty} a_{n-k}x_k-c_n\right)=0, \quad n=0,\pm 1,
  \pm 2 \dots,
\end{equation}
 where \(|r|\neq 1\). With the help of a \(Z\)-transform this can be
 reduced to a boundary value problem of Carleman type~\cite{Ga-Che}.

There are other methods for solving equation of type \eqref{eq:DHol}
without the use of Wiener--Hopf method.
The discrete systems of equations with a Toeplitz operators have effective numerical methods of solutions, some which use
Cauchy-like matrices, which have small numerical
ranks~\cite{Toeplitz_fast_07}.
Also there are some results about perturbation of
Toeplitz operators for example solution of infinite systems like
\(Au+Bu=f\) where \(A\) is Toeplitz operator and \(B\) is in some sense small
\cite{Nowak09}.

Equations of type~\eqref{eq:DHol} naturally appear in applications when
the boundary conditions are discrete and periodic in semi-infinite
configurations. For example, this happens in the Sommerfeld half plane
problem for discrete Helmholz equation \cite{SHARMA17_Somm}. This is also the case
when an semi-infinite discrete array of point scatterers is
considered~\cite{Thompson_14,Albani_19, Mov_Cras_17, Scatterers_15}. Additionally this type of problems is common in crack propagation problems~\cite{Mishuris09}.

\subsection{ Non-uniqueness of factorisation}
\label{subsec:N-U}

This section gives the details of the degree of freedom when constructing
Wiener--Hopf factorisations.  In the scalar case, factors are unique up to
a constant. In other words if there are two such factorisations
\(K(z)=K_+(z)\lambda(z)K_-(z)\) and \(K(z)=P_+(z)\lambda(z)P_-(z)\),
where $\lambda(z)$ is the factor accounting for the index
of the function $K(z)$,
\[K_+=cP_+, \quad \text{and} \quad K_-=c^{-1}P_-,\]
where \(c\) is some non-zero complex constant. This can be seen by
applying analytic continuation to \(K_+/P_+\) and \(P_-/K_-\) and then using
the extended Liouville's theorem~\cite{Kranzer_68}. Here we deal with factorisation of the matrix defined on a simple closed curve $\mathcal L$ in the complex plane. The corresponding result in the \(2\times 2\) matrix case is rather more complicated:

\begin{thm}~\cite{Spit}
\label{thm:Uni}
Let the matrix \(\mathbf{G}( t)\) admit the Wiener--Hopf factorisation~\eqref{eq:F1}. Then:
\[\mathbf{G}=\mathbf{G}_+^{1}\mathbf{\Lambda}\mathbf{G}_-^{1}, \]
  where
  \[\mathbf{G}_+^{1}=\mathbf{G}_+\mathbf{H}, \quad \mathbf{G}_-^{1}=\mathbf{\Lambda}^{-1}\mathbf{H}^{-1}\mathbf{\Lambda}\mathbf{G}_-,\]
  is also a factorisation of \(\mathbf{G}( t)\) where \(\mathbf{H}\) is a constant invertible
  matrix function  if \(\kappa_1=\kappa_2\) and otherwise is of the form:
  \[ \mathbf{H}(t)=\left(
 \begin{array}{cc}
c_1 &  P(z) \\
0  & c_2
 \end{array} \right),\quad z=\frac{t-i}{t+i},
 \]
where \(P(z)\) is a polynomial of degree   \(\kappa_1-\kappa_2\) in $z$.
\end{thm}
This means that there is more  freedom to choose the Wiener--Hopf
factors when \(\kappa_1 \ne \kappa_2\). In the case of higher order matrix functions, similar results can be shown, while the matrix structure has a bit more complicated structure~\cite{Spit_piece-wise}.



\subsection{General class of PDEs that can be reduced
to a Wiener--Hopf}
\label{appendix A}
As an illustration of the concept, we consider a class of partial differential
equations (PDEs) that can be reduced to a Wiener--Hopf equations or a Riemann--Hilbert
boundary value problem~\cite{Ga-Che}*{\S~20.3}. This method is suitable for linear
PDEs of the form:
\begin{equation}
  \label{eq:PDE}
  \sum_{p=0}^A\sum_{q=0}^B h_{pq}(y)\frac{\partial^{p+q}
    u(y,x)}{\partial y^p \partial x^q}=g(y,x),
\end{equation}
where \(h_{pq}\) and \(g\) are given and  \(u\) is to be determined;
note that \(h_{pq}\) has no \(x\) dependence. In the simplest case the
boundary conditions are given on the lines \(y=\)const or half lines
\(y=\)const for \(x>0\) or \(x<0\). Let these \(n\) lines  be positioned at
constants \(y_1\), \dots \(y_n\).  The boundary conditions have the form
\begin{align}
\label{eq:bc}
  \sum_{p=0}^P\sum_{q=0}^Q \left(a_{pqrs}(y)\frac{\partial^{p+q}
    u(y_s+0,x)}{\partial y^p \partial x^q}+b_{pqrs}(y)\frac{\partial^{p+q}
    u(y_s-0,x)}{\partial y^p \partial x^q}\right) =g_{rs}(x), \quad x>0 \\
 \sum_{p=0}^P\sum_{q=0}^Q\left( c_{pqrs}(y)\frac{\partial^{p+q}
    u(y_s+0,x)}{\partial y^p \partial x^q}+d_{pqrs}(y)\frac{\partial^{p+q}
    u(y_s-0,x)}{\partial y^p \partial x^q}\right)=g_{rs}(x), \quad x<0 \\
r=1, \dots m_s, \quad\quad s=1, \dots n.
\end{align}
There could  also be some conditions for \(y\to \infty\) or \(y\to
-\infty\). In order to obtain a Riemann--Hilbert boundary value problem a list of
systematic steps is performed in detail in~\cite{Ga-Che}*{\S~20.3}, and
which is translated into English in~\cite{mythesis}.

We also refer an interested reader to the classic monographs, for
example~\cite{Vekua, Muskh_book,bookWH} for applications of various integral transforms to reduce PDE boundary value problems to Wiener--Hopf or Riemann--Hilbert form, amongst others.

\subsection{Applications}
\label{sec:app}

The wide applicability of the Wiener--Hopf methods has ensured its
continuous development from its inception
to modern day.
Traditional applications areas are acoustics, aeroacoustics, water
waves, Levi processes, signal processing,
electromagnetism~\cite{Application_DK, sto-W-H-eng,vorovich1979,
  Kopev08, daniele2014wiener,finance}. Various static and dynamic
problems of fracture mechanics~\cite{Piccolroaz07,Pricc09,Antipov99,
  Willis95,Slepyan15}, lattices~\cite{Ayzenberg14, Slepyan02,
  Brun14,Maurya19, Sharma20, Shanin20,GMT} and recently in nanophotonics,
matamaterials and biomechanics are studied using this
technique~\cite{Mishuris09,nano_Wiener_19, Albani_meta,Ant_crack,GT}.

The Wiener--Hopf method is extensively used for canonical scattering
problems both in acoustics~\cite{bookWH} and
electromagnetism~\cite{daniele2014wiener}. In both of these
applications the governing equation is Helmholtz' equation.
There are fewer methods available to tackle Helmholtz' equation compared
to Laplace's equation. For the Laplace equation, other techniques from complex
analysis such as conformal mapping are applicable \cite{Crowdy_con}.
Although the Helmholtz equation is self-adjoint
addition of boundary conditions, including the Sommerfield radiation
condition, means that the boundary value problem is no longer
self-adjoint. This causes problems in demonstrating completeness of series
expansions \cite{FokSpen12}. The Wiener--Hopf method
on the other hand very naturally incorporates the radiation condition
and complicated edge conditions.  Simple models in aeroacoustics can also be reduced
to acoustic problems via Lighthill's analogy and the reciprocal theorem~\cite{Lighthill}.

Water waves is an area where the Wiener--Hopf techniques is accepted as
one of key analytic
tools~\cite{Handbook_water}. The ice cover on water can be modelled as
an elastic material for
 which flexural motions are able to be described by thin-plate or Timoshenko-Mindlin
 theory~\cite{Craster_ice, Smith_20}.
Eigenfunction-matching method is also a
commonly used method and in some cases this has been shown to be
equivalent to the Wiener--Hopf method~\cite{Water_equi}.
Recently, the problem of waves generated in a fluid and an ice sheet by a pressure region moving on the free surface of the fluid, along the edge of the semi-infinite ice sheet, has been solved using the Wiener--Hopf technique~\cite{ice_free_surface}.

One of the new applications of the Wiener--Hopf technique has been
nanophotonics. Graphene and other two-dimensional (2D) materials may  sustain evanescent, fine-scale electromagnetic waves that are tightly
confined to the boundary, called
plasmonpolariton~\cite{nano_Wiener_19}. This is also closely related to
edge magnetoplasmons, the theory of which was also derived using the
Wiener--Hopf method~\cite{Mag_edge_88} and the more realistic extension of
which lead to more complicated Wiener--Hopf
systems~\cite{magneto_18}. Recently, composites and Metamaterials have also been investigated using the Wiener--Hopf method~\cite{Albani_meta,Thompson_14,Bankov2008, Gower_19}.

This is by no means an exhaustive list of applications; unfortunately space constrains us in this short survey. We now list a number of approaches to particular
constructive exact or approximate factorisation methods in the sections below.

\section{A list of constructive procedures}
\label{constructive}

By constructive factorisation we mean
that there is an algorithmic way to obtain the
factorisation  \eqref{eq:F1} or
\eqref{WH_fact}. Recall, for scalar functions, factorisation is always possible
via integrals of Cauchy type \eqref{eq:R5} and \eqref{eq:R6}.
Currently, there are no constructive factorisation
procedures for matrix functions in general, and only very specific classes can be exactly solved.
What complicates the development of an algorithmic approach is
that for a given a matrix function it is difficult to determine if any known
constructive procedure apply. Pre- and post-multiplication can sometimes transform the matrix to a
known form, but this is undertaken mainly using \emph{ad hoc} techniques. A
systematic treatment of some transformations is given in \cite{Ehr_Spec}.

The main classes of matrix kernels which have a constructive factorisation are reviewed below.
 They are the rational matrix functions, commutative matrices particular matrices
 with exponential factors, as well
 as several miscellaneous classes.
We also refer the reader to several other important classes of piece-wise constant
matrix functions not included in the review,
see~\cite{Spit_piece-wise, Duduchava79,Rogosin_21}, the Meister and Speck
matrix as discussed in \cite{Abrahams_02}
and for applications see~\cite{Antipov_09, Antipov_18}. We also do not
consider the case where there are singularities on the curve
\({\mathcal L}\)~\cite{Mikhlin86}; for more information the reader is
refereed to~\cite{Gohberg74, Prossdorf91}.

\subsection{Rational matrix function}
\label{rational}

This is the simplest class of matrix functions for which an exact constructive
factorisation is known. This is because there are only a finite number
of isolated singularities which need to be considered. A procedure to construct
the Wiener--Hopf factorisation is known in this case and is briefly described
immediately below, followed by the standard ``pole removal'' technique.

\subsubsection{Factorisation of rational matrix functions}
\label{subsec:rational}

This subsection describes an algorithm to
factorise  a rational matrix function ~\cite{Spit}*{Sec 1.2}. Given a matrix function
\(M\) it can be expressed as \(P/q\) where \(P\) is a polynomial
matrix function and \(q\) is a scalar polynomial.  It is enough to
factor them separately. The scalar function \(q\) can be factorised easily, see Section~\ref{sec:R-H}. We will assume that  \(\det P\) has no roots
on the real axis. What prevents a polynomial matrix \(P\) being a plus factor
is the presence of zeros of \(\det P\) in the upper half-plane.
(Recall that, for factorisation, a plus (minus) factor requires the matrix \textbf{and} its inverse to be analytic in the upper (lower) half plane.)
A way of eliminating one zero at a time is outlined in the following lemma.

\begin{lem}
Let \(P\) be a polynomial matrix function such that \(\det\, P\) has no roots
on the real axis. If \(\det\, P\) has \(n\) roots in the upper
half-plane then we can construct a polynomial \(L=PR\) where \(\det\, L\) has
\(n-1\) roots in the upper half-plane and the matrix \(R\) has the form (blanks stand for zeroes):
\begin{equation}
  \left(
 \begin{array}{cccccc}
1 & & &\frac{-c_{1}}{z-z_0} &&\\
 & 1& & \vdots&&\\
 && \ddots &\frac{-c_{k-1}}{z-z_0} &&\\
 &&& \frac{1}{z-z_0} &&\\
 &&&  &\ddots&\\
 & & && &1
\end{array}\right),
\end{equation}
where $z=z_0$ is one of the aforementioned poles, while $c_k$ are computed constructively and uniquely.
\end{lem}

We can repeat this procedure, eliminating all the roots of \(\det\, P\)
in the upper half-plane and hence construct the plus factor. Note that the
matrix \(R\) is invertible.  The minus factor can be found either by pre-multiplying \(P\) by the inverse of the plus factor, or by taking the inverse of the product of all of the \(R\) matrices. This completes the factorisation.

\subsubsection{ ``Pole Removal'' for rational matrix functions}
\label{subsec:pole}

In applications the above rational factorisation algorithm is
rarely used; instead ``pole removal'' or ``singularity matching'' is
employed~\cite{Abraha_all_pade, daniele2014wiener} and is recapped below.

Suppose that  \(A(z)\) is rational (every entry of the matrix function
is a rational function), then we can solve~\eqref{eq:WH2}
without multiplicative factorisation. Perform an additive split \(A(z)\Phi_+(z)=(A(
z)\Phi_+(z))^-+(A( z)\Phi_+(z))^+\) and \(C(z)=C^-(z)+C^+(z)\), where the first term is analytic
in ${\mathcal H}^{-}$ and second in ${\mathcal H}^{+}$.
Since \(A(z)\) is rational we can write \((A(
z)\Phi_+(z))^-= \sum_{i=0}^n \frac{A_i}{z-z_i}\) where \(z_i\) are the poles in the
upper half-plane, and \(A_i\) are constants (residues at \(z_i\)). Thus,
\[A( z)\Phi_+(z)-\sum_{i=0}^n
\frac{A_i}{z-z_i}+C^+(z)= -\Psi_-(z)-\sum_{i=0}^n
\frac{A_i}{z-z_i} -C^-(z)\]
separates terms on the left on the left and right hand sides into functions
analytic in the upper and lower half planes respectively.  To find \(A_i\)
we need to solve a linear system of equations (by setting \(z=z_i\)).

The pole removal method has an advantage over the rational
factorisation procedure. It enables one to remove undesired poles without changing
the rest of the equations. This is one of the reason that it can not
only be applied to rational Wiener--Hopf kernels but to more general
kernels which have one or more rational parts. Pole removal is closely
linked to the generalised Liouville's theorem (which allows the unknown
function to have poles). It can also be naturally extended to
meromorphic functions by allowing an infinite sum, which would lead to
an infinite linear system that would need to be truncated in order to
be solved.
The additional advantage of ``pole removal'' is that the
multiplicative factorisation is reduced to an additive splitting. It is also worth noting that many alternative formulations, such as Fredholm factorisation~\cite{daniele2014wiener, Daniele_Lombardi_11, Daniele_Lombardi_16},  approximate techniques (see Section~\ref{sec:approx}) and numerical methods~\cite{trogdon_riemannhilbert_2016}, are also based on this
principle or replacing a multiplicative factorisation by additive splittings.

\begin{rem}
Note that the instability issue discussed
in Section~\ref{sec:approx} still applies here. To recognise the challenge it is enough to refer to the fact that, generally speaking, all poles or zeros are computed approximately,
while the system of linear equations mentioned above (which have to be solved to compute the constants) can be ill-conditioned. In case of a stable set of the partial indices the system is not ill-conditioned, otherwise, at least for some of the poles, this would be expected~\cite{Adukov_20}.
\end{rem}

\subsection{Analytic and meromorphic matrix functions}
\label{analytic}

In \cite{Adu93} the factorisation problem for the meromorphic matrix function was solved. The corresponding algorithm is based on the reduction of the problem to a finite number of systems of linear algebraic equations, whose coefficients' matrices are block Toeplitz matrices:
\begin{equation}
\label{an1} T_k = \left(\begin{array}{cccc} C_{k} & C_{k-1} &
\cdots & C_{-2 \kappa} \\
C_{k+1} & C_{k} &
\cdots & C_{-2 \kappa+1} \\
\cdots & \cdots &
\cdots & \cdots \\
C_{0} & C_{-1} & \cdots & C_{-2 \kappa - k}
\end{array}
\right),
\end{equation}
with $C_j$ being the power moments of the inverse $G^{-1}(t)$ to the given matrix  with respect to the contour ${\mathcal L}$
\begin{equation}
\label{an2}
C_j = \frac{1}{2\pi i} \int\limits_{\mathcal L} t^{-j-1} G^{-1}(t) dt.
\end{equation}
The partial indices of the left and right factorisation are
represented in terms of ranks of the above Toeplitz
matrices. An analogous method is applied in \cite{Adu99} for the case of
analytic matrix functions which, being a special case of the latter,
need less computation. A part of this algorithm related to factorisation of the scalar polynomial is
implemented in the Maple computer system as module
\verb!PolynomialFactorization ! (see its description and illustrative examples in \cite{Adu18}).

\subsection{Triangular matrix functions}
\label{triangle}

Upper or lower triangular matrix functions (with factorisable diagonal
elements) can  sometimes be reduced to a set of scalar
 equations. If the equations can be solved in turn and the previous
 solution used to make the next equation scalar, then the system
 decouples. However, even in the case of triangular matrix functions there
 may be complications which mean that the system does not decouple in this way. This happens when the partial indices are
 non-zero, for which case we describe a general procedure below. It relies on the
 idea of a canonical matrix, introduced in Section~\ref{sec:pre}\ref{subsec:can_mat}
 and this was one the first attempts to realize this construction explicitly
 (i.e.\ to determine canonical matrix factorisation for (\ref{R0})).
The classical results by Chebotarev  \cite{Che56} applies to $2\times 2$ triangular matrix functions with factorisable diagonal entries.

Let
$$
G(t) = \left(\begin{array}{cc}
\zeta_1(t) & 0 \\
a(t) & \zeta_2(t)
\end{array}
\right), \quad t\in {\mathbb T} \quad \text{(unit circle)}.
$$
Denote $\kappa_j = {\mathrm{ind}}_{\Gamma} \zeta_j(t)$  and let $x^{\pm}_j(z)$ be canonical functions for the
scalar homogeneous Riemann--Hilbert boundary value problems
with coefficients $\zeta_j(t)$, respectively (see \cite{Gak77}), $j = 1, 2$. Then the piece-wise analytic matrix
$$
X^{\pm}(t) = \left(\begin{array}{cc}
x_1^{\pm}(t) & 0 \\
x_2^{\pm}(t) \phi^{\pm}(t) & x_2^{\pm}(t)
\end{array}
\right)
$$
with
$$
\phi^{\pm}(z) = \frac{1}{2\pi
i} \int\limits_{\mathbb T} \frac{a(\tau) x_1^{-}(\tau)
d\tau}{x_2^{+}(\tau) (\tau - z)}, \; z\in D^{\pm},
$$
satisfies the following boundary condition:
$
X^{+}(t) = G(t) X^{-}(t).
$
Let $\mu \geq 1$ be the order of the function $\phi^{-}(t)$ at infinity. The orders of non-zero elements of the matrix $X^{-}(z)$ at infinity
can be characterized by the following table
$$
\left(\begin{array}{cc}
\kappa_1 &  - \\
\kappa_2 + \mu & \kappa_2
\end{array}
\right).
$$

If $\kappa_1 \leq \kappa_2 + \mu$, then the matrix $X^{-}(z)$ has
the normal form at infinity, i.e. it is the canonical matrix. Thus,
partial indices of the matrix  $G(t)$ are equal $(\kappa_1,
\kappa_2)$. In the case $\kappa_1 > \kappa_2 + \mu$ Chebotarev
proposed to use an expansion of $\frac{1}{\phi^{-}(z)}$ as a continued fraction
$$
\frac{1}{\phi^{-}(z)} = q^{\gamma_0}(z) + \frac{1}{q^{\gamma_1}(z) + \frac{1}{q^{\gamma_2}(z) + \ldots}},
$$
where $q^{\gamma_i}(z)$ are polynomials of orders $\gamma_i$, respectively ($\gamma_0 = \mu$). Using these polynomials in the scheme of the elementary transformations it was shown in \cite{Che56} that in a finite number of steps the ``minus''-matrix can be rebuilt to the normal form at infinity and thus the partial indices will be found and factors in \eqref{eq:F1} will be determined.
In \cite{FGK95} this approach was applied to obtain an explicit solution of
the factorisation problem for a case of the Khrapkov-Daniele matrix
functions. In \cite{Primachuk_20} this approach was applied to solve \(\mathbb{R}\)-linear conjugation problem.

In \cite{PriRog18} Chebotarev's approach was realized for triangular
matrix function of arbitrary order with factorisable diagonal
entries.  The basis for the inductive steps is the following statement:
{\it let ${\mathcal L}$ be a simple smooth closed contour $0\in
  {\mathcal L}$, and let $B(t)$, for $t\in {\mathcal L},$ be a
non-singular H\"older continuous square matrix-function of order $n$ having the following form:
\begin{equation}
\label{PriRog1} B(t) = \left( \begin{array}{cc} A(t) & {\bf 0}\\
b_{1}(t) \ldots b_{n-1}(t) & c(t)
\end{array}
\right), \;\;\; {\bf 0} = \left(\begin{array}{c} 0 \\ \vdots \\
0\end{array}\right).
\end{equation}
Suppose that the non-singular square matrix-function $A(t)$ of the order $n-1$ admits the factorisation
$$
A(t) = A^{+}(t)\; \Lambda(t)\; A^{-}(t),
$$
where $\Lambda(t) = {\mathrm{diag}}\, \left\{t^{{\kappa}_1}, \ldots,
t^{{\kappa}_{n-1}}\right\}$. Then the matrix-function $B(t)$
possesses a factorisation if the following matrix does:
\begin{equation}
\label{PriRog2}
\left( \begin{array}{cc} \Lambda(t) & {\bf 0}\\
({\bf b}(t)|{\bf Y}_1^{-}(t))  \ldots ({\bf b}(t)|{\bf
Y}_{n-1}^{-}(t)) & c(t)
\end{array}
\right).
\end{equation}
Here ${\bf b}(t) = (b_1(t),\ldots,b_{n-1}(t))$ is the row of the first
${n-1}$ entries of the lowest row of $B(t)$,  ${\bf Y}_j^{-}(t) =
(y^{-}_{1j}(t),\ldots, y^{-}_{(n-1)j}(t))^{T}$ is the  $j$th column of the matrix-function
$Y^{-}(t) = (X^{-}(t))^{-1}$, and
$$
({\bf b}(t)|{\bf Y}_j^{-}(t)) = \sum\limits_{k=1}^{n-1} b_k(t)
y^{-}_{kj}(t).
$$
}
This process can be repeated \(n\) times. 

Note that for a matrix with a specific structure, the general stability criterion for the partial indices can be extended (in other words, within a specific sub-algebra, there are cases when the perturbation preserving the structure remains stable, while there exists a general (small) perturbation of other structure that changes the set of the partial indices. For the 2x2 triangular matrices such analysis has been conducted in~\cite{Adukov_20}. 

\subsection{Functionally commutative matrix functions}
\label{FunCom}

There has been significant progress in the case of the matrices
possessing commutative factorisations~\cite{Spit_book}. One of
the reasons is that one can apply some of the
techniques that have been developed in the scalar
setting. We present the
corresponding result in the case of H\"older continuous coefficients
$G$ and $g$ in~\eqref{eq:R1}.

 In \cite{Gak52}, the following question was considered: whether there
 exist a classes of matrix coefficients
$G(t)$ such
that formula \eqref{eq:R3} still gives the bounded solution of the vector-matrix problem \eqref{eq:R1} (i.e. for
$n  > 1$)? Such a question is related to the purely algebraic question of the fulfillment of the identity
\begin{equation}
\label{FunCom1}
e^{A(t)}\cdot e^{B(t)} = e^{A(t) + B(t)},
\end{equation}
for certain matrices $A(t)$ and $ B(t)$ connected to the coefficient
$G(t)$. It was shown in \cite{Gak52} that
identity (\ref{FunCom1}) holds (and thus Gakhov's formula \eqref{eq:R3} can be used to solve problem \eqref{eq:R1} in the matrix case), whenever the matrix function $G(t)$ is functionally commutative. The latter
condition means
\begin{equation}
\label{FunCom2}
G(t)G(\tau) =G(\tau)G(t), \; \forall \; t, \tau\in {\mathcal L}.
\end{equation}
The question of solvability of the vector-matrix problem \eqref{eq:R1} in the case of a functionally commutative matrix
$G(t)$ is solved in \cite{Gak52}. In \cite{Che56a} the following necessary and sufficient condition for a matrix to be functionally commutative is reported: {\it  the matrix $G(t)$ is functionally commutative if and only iff it can be represented in the form}
\begin{equation}
\label{FunCom3}
G(t) = \varphi_1(t) G_1 + \varphi_2(t) G_2 + \ldots + \varphi_m(t) G_m,
\end{equation}
{\it where $\varphi_1(t), \varphi_2(t), \ldots, \varphi_m(t)$ are
  linear scalar independent functions, and $G_1, G_2, \ldots, G_m$ are
  constant linear independent mutually commutative matrices.} The
structure of lower-dimensional (up to $n = 4$) functionally commutative
matrix functions was described in  \cite{Che56a} too. In
\cite{RogDub17} a solution to the Riemann--Hilbert boundary value problem
\eqref{eq:R1} and the corresponding factorisation problem is
illustrated by examples (see also \cite{Adu13}).

\subsection{Commutative factorisation}
\label{CommFact}

Close to the above case is so called commutative factorisation. Let us
formulate it in the Wiener--Hopf setting. Let the matrix $K(z)$ be defined on
a strip ${\mathcal H}= \{ z \in \mathbb{C}: a < \mathrm{Im}\, z < b\}$  around the real line. The
{\it commutative factorisation} of \(K(z)\) are factors \(K_-(z)\)
and \(K_+(z)\)  satisfying \eqref{WH_fact} and the factors $K_{+}$,
$K_{-}$ commute in \({\mathcal H}\). The idea of commutative factorisation stems from the work of Heins \cite{Hei50}.
For criteria  when
commutative factorisation is possible the reader is referred to \cite{comm}. Several special
cases of matrices admitting commutative factorisations are described below.

\subsubsection{Khrapkov-Daniele matrices}
\label{KhrDan}

Matrices of this kind represent one of the  simplest classes of
matrices with commutative factorisation. In $2\times 2$ case these
matrices were discovered from certain diffraction problems
\cite{Khr71a,Khr71b,Dan78}. They appear in the study of diffraction by
geometries which had a right angle such as wedges or
gratings (a grating can be thought of as a perpendicular strip in a
waveguide)~\cite{gratting} as well as other
configurations~\cite{Application_DK,Khr01,VeiAbr07,DK_plas,Antipov_DK}.
The Khrapkov-Daniele matrices have the following form
\begin{equation}
\label{KhrDan1}
K(\alpha) = k_0(\alpha) I + k_1(\alpha) J(\alpha),
\end{equation}
where $J(\alpha)$ is a polynomial matrix-function with the property
\begin{equation}
\label{KhrDan2}
J^2(\alpha) = \Delta^2(\alpha) I,
\end{equation}
with $k_0, k_1$ being arbitrary functions analytic in \({\mathcal H}\) and having algebraic growth at
infinity and $\Delta^2$ is a polynomial in $\alpha$,
$$
\Delta^2(\alpha) = O\left(|\alpha|^p\right), \; \alpha \rightarrow \infty,
$$
with $p \leq 2$.

Commutative factorisation under the above conditions has the form (see e.g. \cite{VeiAbr07})
\begin{equation}
\label{KhrDan3}
K^{\pm}(\alpha) = r^{\pm}(\alpha) \exp\left\{\theta^{\pm}(\alpha) J\right\},
\end{equation}
where scalar functions $r^{\pm}, \theta^{\pm}$ satisfy the relations
\begin{equation}
\label{KhrDan4}
r^{+}(\alpha) r^{-}(\alpha) =r(\alpha) := \sqrt{k_0^2(\alpha) - \Delta^2(\alpha) k_1^2(\alpha)},
\end{equation}
\begin{equation}
\label{KhrDan5}
\theta^{+}(\alpha) + \theta^{-}(\alpha) = \theta(\alpha) := \frac{1}{\Delta(\alpha)} \tanh^{-1} \left(\frac{k_1(\alpha)}{k_0(\alpha)} \Delta(\alpha)\right).
\end{equation}

An interesting approximate technique, based on the Khrapkov-Daniele matrix form, which allows it to be extended to certain non-commutative matrices, has also been developed \cite{Elastic_half}.

\subsubsection{Jones class}
\label{Jones}

Jones \cite{Jon84} established the commutative factorisation of the following $n\times n$ class of matrices
\begin{equation}
\label{Jones1}
C(\alpha) = \sum\limits_{m=1}^{n} a_m(\alpha) E^m(\alpha),
\end{equation}
where $a_1(\alpha), a_2(\alpha), \ldots, a_n(\alpha)$ are functions analytic in the strip \({\mathcal H}\), $E(\alpha)$ is an entire matrix-function, $E^n = q^n(\alpha) I$. It is supposed also that $\det\, C(\alpha)\not= 0$ and ${\mathrm{tr}}\, E^r = 0, r =1,2,\ldots,n-1.$

Under these conditions it was shown that $C(\alpha)$ possesses the commutative factorisation
$$
C(\alpha) = C^{+}(\alpha) C^{-}(\alpha) = C^{-}(\alpha) C^{+}(\alpha),
$$
where
\begin{equation}
\label{Jones2}
C^{\pm}(\alpha) = \left[(\det\, C(\alpha))^{\pm}\right]^{1/n} \sum\limits_{p=0}^{n-1} \gamma_p^{\pm}(\alpha) E^p(\alpha),
\end{equation}
\begin{equation}
\label{Jones3}
\gamma_p^{\pm} = \sum\limits_{m=0}^{n-1} \frac{\omega^{m p}}{n q^p} \exp \left(\sum\limits_{r=1}^{n-1} q^{n - r} \omega^{r m} \delta_{n - r}^{\pm}\right),
\end{equation}
\begin{equation}
\label{Jones4}
\delta_s^{+} + \delta_s^{-} = \delta_s :=   \sum\limits_{p=0}^{n-1} \frac{\omega^{s p}}{n q^s} \ln \left[\sum\limits_{l=1}^{n} q^{n - l} \omega^{l p} a_{n - l}\right].
\end{equation}
Here $\omega = e^{(2\pi i)/n}$ and the functions $\delta_s^{+}, \delta_s^{-}$ are analytic in the corresponding half-planes.

In \cite{VeiAbr07} the Wiener--Hopf factorisation is constructed for
matrices of the Jones class under assumption that $E$ is an entire matrix-function with polynomial elements having distinct eigenvalues $\lambda_1, \lambda_2, \ldots, \lambda_n$ such that the following condition holds
\begin{equation}
\label{Jones5}
(E^{n_1} - \mu_1 I) (E^{n_2} - \mu_2 I) \ldots (E^{n_p} - \mu_p I) = 0,
\end{equation}
where $\mu_j$ are polynomials and $n_1 + n_2 + \ldots n_p = n$.

\subsubsection{Moiseev class}
\label{Moiseev}

In \cite{Moi89,Moi93} the following class of matrices (similar to the Jones class but satisfying different assumptions) was considered
\begin{equation}
\label{Moiseev1}
G(\alpha) = \sum\limits_{m=0}^{n-1} a_m(\alpha) E^m(\alpha),
\end{equation}
where $a_m(\alpha)$ are scalar functions and $E(\alpha)$ is a polynomial matrix. In the simplest
case of matrix $E$ having distinct eigenvalues almost everywhere, $E$ can be
decomposed as
\begin{equation}
\label{Moiseev2}
E =T \Lambda T^{-1},
\end{equation}
where $T$ is the matrix of the eigenvectors and $\Lambda$ is a diagonal matrix composed
of the eigenvalues. Both $T$ and $\Lambda$ are algebraic matrices, and in \cite{Moi89}  the Riemann surface ${\mathcal R}$ was
introduced on which $T$ and
$\Lambda$ are single valued. Further,
the matrix factorisation problem reduced to a scalar Riemann--Hilbert boundary value problem on ${\mathcal R}$. This problem can be solved in terms of Abelian integrals
with the help of Jacobi's inversion problem as stated by Zverovich in \cite{Zve71}. So, the solution to the
problem of factorisation of (\ref{Moiseev1}) is known at least formally, and it possibly can
be used for practical purposes. In some special cases this scheme was realized in a series of articles by Antipov and
Silvestrov, and the solution to the factorisation problem
(as well as to the corresponding problems arising in application) was constructed in an explicit form in terms
of special functions (see \cite{AntSil02, Ant14} and references
therein). The idea of reducing a matrix factorisation to a scalar one
on a Riemann surface is also developed in \cite{Camara_sur}.

In \cite{Shanin_13} a simplified method for treating matrices
of  Moiseev's  type is proposed. An additional restriction is that
$a_m(\alpha)$ are algebraic functions. The factorisation problem is
transformed by Hurd's procedure (see, e.g. \cite{Hur76,Hur87,daniele2014wiener})  into a
Riemann--Hilbert boundary value problem on a set of cuts. The method relies on formulation of ordinary differential
equation with an unknown coefficient. It is shown that the proposed procedure
in the Khrapkov case is equivalent to the standard solution method.

\subsection{Factorisation of matrices with exponential terms in the entries}
\label{exponent}


\subsubsection{Integral equation on a finite interval}
\label{int}

Consider the  system of the convolution equations on a finite interval
\begin{equation}
\label{int1}
\varphi(t) + \int\limits_{0}^{a} k_0(t - s) \varphi(s) ds = f(t), \; 0 < t < a,
\end{equation}
note the similarity and differences with the basic integral
Wiener--Hopf equation \eqref{eq:WH1}. The application of the Fourier transform leads to the necessity to factorize the $2n\times 2n$ matrix function of type
\begin{equation}
\label{int2}
A(\lambda) = \left(\begin{array}{cc} - e^{-i a z} K(z) & - I + K(z)\\
I + K(z) & - e^{i a z} K(z)\end{array}\right),
\end{equation}
where the $n\times n$ matrix $ K(z)$ is the Fourier transform of any extension $k$ onto the whole real line ${\mathbb R}$ of the kernel $k_0$.

The factorisation of matrix (\ref{int2}) is studied in \cite{FGK00} under condition that the operator
\begin{equation}
\label{int3}
\left(B \varphi\right)(t) := \varphi(t) + \int\limits_{0}^{a} k_0(t - s) \varphi(s) ds,
\end{equation}
is invertible in $L_{n\times 1}(0, a)$ and it is a-priori known that partial indices of the matrix $A(\lambda)$ are equal to zero.
In this case the corresponding factors are found explicitly.

\subsubsection{A system of integral equations with an exponential kernel}
\label{system}

In \cite{Abr_exp} the following system of two convolution-type equations on the real semi-axes is considered:
\begin{equation}
\label{system1}
u(x) = \int\limits_{0}^{\infty} k(x - t) u(t) dt + f(x),\; x\in {\mathbb R}_{+},
\end{equation}
where the kernel $k$ is the following matrix
\begin{equation}
\label{system2}
k_b(x) = \left(\begin{array}{cc} e^{-|x|} & e^{-|x - b|} \\ e^{-|x + b|} & e^{-|x|} \end{array}\right).
\end{equation}
Fourier transformation of this system leads to a functional equation with matrix kernel
\begin{equation}
K_b(\alpha) = \left(\begin{array}{cc} 2 \lambda - 1 - \alpha^2 & 2 \lambda e^{i b \alpha} \\ 2 \lambda e^{- i b \alpha} & 2 \lambda - 1 - \alpha^2 \end{array}\right).
\end{equation}
Factorisation of matrices of such type, as well as some other matrices involving exponential elements, is given in~\cite{Abr_exp}.

\subsubsection{AKM-approach}
\label{AKM}

The name of this approach is simply an abbreviation of the authors' names
in the paper \cite{Akto92} which considers factorisation of $2\times 2$ non-rational matrix functions of the form
\begin{equation}
\label{AKM1}
G(t; k) = \left(\begin{array}{cc} T(t) & - R(t) e^{2i kt} \\ - L(t) e^{-2i kt} & T(t) \end{array}\right), \;\;\; t\in {\mathbb R},
\end{equation}
for an arbitrary real parameter $k$. Matrices of this type arise as (modified) scattering matrices  for the 1-D Schr\"odinger equation and some related Schr\"odinger-type equations. For the above discussion the so called scattering matrix
\begin{equation}
\label{AKM1s}
S(t) =  \left(\begin{array}{cc} T(t) &  R(t) \\ L(t) & T(t) \end{array}\right),
\end{equation}
plays an important role. Here $T$ is called the transmission coefficient, and $R, L$ are reflection coefficients from the right and from the left, respectively.

It is supposed that the following conditions are satisfied
\begin{enumerate}
\item  $T(z)\not= 0$, $z\in \overline{\Pi^{+}}\setminus\{0\} = \{z\in
  {\mathbb C}: {\mathrm{Im}}\, z \geq 0, z\not= 0\}$, is meromorphic
  in $\Pi^{+}$ with continuous boundary values on the extended real
  line, either $T(0)\not=0$ or $T(t)$ is vanishing linearly at $t=0$,
  and $T(\infty) = 1$.
\item  $R(z)$ and $L(z)$ are meromorphic on $\Pi^{+}$ with continuous
  boundary values on the extended real line and vanishes as
  $z\rightarrow \infty$ in $\overline{\Pi^{+}}$.
\item $G(t; k)^{-1} = Q G(-t; k) Q$ for all $t \in {\mathbb R}$, where
  $Q = \left(\begin{array}{cc} 0 & 1 \\ 1 & 0\end{array}\right)$.
\item $G(t; k)$, as the function of $t\in {\mathbb R}$, belongs to a suitable Banach algebra of $2\times 2$ matrix functions within which factorisation is possible. This may be the Wiener algebra ${\mathcal W} = {\mathcal W}^{2\times 2}$ or algebra of functions $f(t)$ such that $f^{\ast}(\xi) = f\left(i \frac{1 + \xi}{1 - \xi}\right)\in H_{\alpha}({\mathbb T}), \alpha\in (0, 1)$ (which we denote ${\mathcal H}_{\alpha}^{2\times 2}$).
\end{enumerate}

First the following statement was proved: {\it let us consider
\begin{equation}
\label{AKM2}
W(t)= \left(\begin{array}{cc} 1 & q(t) \\ - \overline{q(t)} & 1\end{array}\right),\;\;\; t\in {\mathbb R},
\end{equation}
where $q(\infty) = 0$, and either $q\in {\mathbb W}$ or $q\in H_{\alpha}({\mathbb T}), \alpha\in (0, 1)$. Then $W(t)$ has a unique (right) canonical factorisation
$$
W(t) = W_{-}(t) W_{+}(t),\;\;\; t\in {\mathbb R},
$$
where either $W_{\pm}\in {\mathcal W}^{2\times 2}$ or $W_{\pm}\in {\mathcal H}_{\alpha}^{2\times 2}$, and $W_{\pm}(\infty) = E_2$.}

It leads to the following result: {\it let
 matrix (\ref{AKM1}) be a unitary matrix for all $t\in {\mathbb R}$
 satisfying the following conditions:

 \begin{enumerate}
 \item $T(t)\not= 0$ for all $t\in {\mathbb R}$,
 \item $T(t)$ can be continued to a meromorphic function on $\Pi^{+}$ with continuous boundary values on $\overline{\mathbb R}$, and $T(\infty) = 1$.
 \item $T(t)$, $R(t)$ and $L(t)$ belong to either ${\mathbb W}$ or to $H_{\alpha}, 0 < \alpha < 1.$
 \end{enumerate}
Then $G(t; k)$ has a (right) factorisation with equal partial indices $\kappa_1 = \kappa_2 = \kappa$, where $\kappa$ is the number of zeros minus the number of poles of $T(z)$ in $\Pi^{+}$.

The factorisation has the form
$$
G(t; k) = W_{-}(t,k) T_{-}(t) {\mathrm{diag}}\, \left\{\left(\frac{t - i}{t + i}\right)^{\kappa}, \left(\frac{t - i}{t + i}\right)^{\kappa}\right\} T_{+}(t) W_{+}(t,k),
$$
where $W_{-}, W_{+}$ are factors of the canonical decomposition of the matrix $\frac{G(t;k)}{T(t)}$ and $T_{-}, T_{+}$ are product factors of the scalar function $T(t)$:
$$
T(t) = T_{-}(t) \left(\frac{t - i}{t + i}\right)^{\kappa} T_{+}(t).
$$
}
Some approximate methods for matrices with exponential entries are considered in Section~\ref{sec:approx}\ref{sec:Neg}.

\subsection{Miscellaneous Classes}

\subsubsection{(EF)-Algorithm}
\label{EF}

An algorithm based on the successive solution of the scalar boundary value problems was proposed by Feldman, Gohberg and Krupnik \cite{FGK94}.
They construct the solution of the factorisation problem \eqref{eq:F1}
when the curve ${\mathcal L}$ is either the unit circle ${\mathbb T}$
or the real line ${\mathbb R}$.

In the case ${\mathcal L} = {\mathbb T}$
the $2\times 2$ matrix investigated is a non-singular matrix of the following type
\begin{equation}
\label{EF1}
G(t) = \left(\begin{array}{cc} 1 & b(t) \\ c(t) &
                                                  d(t) \end{array}\right), \quad
                                              t \in \mathbb T,
\end{equation}

satisfying the following conditions:
\begin{enumerate}
\item $b, c, d$ belong to the Wiener algebra ${\mathcal W}({\mathbb
    T})$ of absolutely converging Fourier series on ${\mathbb T}$;
\item $b(t) = \frac{p(t)}{q(t)}$, where $p\in {\mathcal W}_{+}({\mathbb T})$ (i.e. is a Fourier series with vanishing coefficients of negative indices), and $q$ is a polynomial $q(t) = \prod\limits_{j=1}^{k} (t - a_j)^{m_j}$ with distinct zeros $a_j$ lying in the unit disc, $|a_j| < 1$.
\end{enumerate}

In the case ${\mathcal L} = {\mathbb R}$ the Wiener algebra ${\mathcal W}({\mathbb R})$ (and its subalgebras ${\mathcal W}_{+}({\mathbb R})$, ${\mathcal W}_{-}({\mathbb R})$) of the functions $\phi(\lambda), -\infty < \lambda < +\infty,$ are defined in the standard way
\begin{equation}
\label{EF2}
{\mathcal W}({\mathbb R}) \ni \phi(\lambda) = c + \int\limits_{-\infty}^{+\infty} \varphi(\tau) e^{i \lambda \tau} d\tau,
\end{equation}
$$
{\mathcal W}_{+}({\mathbb R}) \ni \phi(\lambda) = c + \int\limits_{0}^{+\infty} \varphi(\tau) e^{i \lambda \tau} d\tau,
$$
$$
{\mathcal W}_{-}({\mathbb R}) \ni \phi(\lambda) = c + \int\limits_{-\infty}^{0} \varphi(\tau) e^{i \lambda \tau} d\tau,
$$
where $c$ is an arbitrary complex constant and $\varphi\in L_1$ on its domain of integration.

The matrix under consideration is the following:
\begin{equation}
\label{EF3}
G(t) = \left(\begin{array}{cc} a(t) & b(t) \\ c(t) &
                                                     1 \end{array}\right),
                                                \quad t \in \mathbb R,
\end{equation}

with the corresponding conditions
\begin{enumerate}
\item $a, b, c$ belong to the Wiener algebra ${\mathcal W}({\mathbb
    R})$;
 \item $b(t) = \frac{p(t)}{q(t)}$, where $p\in {\mathcal W}_{+}({\mathbb R})$, and $q$ is a rational function
$q(t) = \prod\limits_{j=1}^{k} \left(\frac{t - a_j}{t + i \gamma_j}\right)^{m_j}$ with distinct zeros $a_j$ laying in the upper half-plane, i.e.  ${\mathrm{Im}}\, a_j > 0$, and $\gamma_j > 0$.
\end{enumerate}

In \cite{FGK04} this algorithm is applied to the solution of two problems: a) the factorisation of matrix  (\ref{EF1}) with ratio $f(t) = \frac{p(t)}{q(t)}$ being meromorphically extended either inside the unit disc or outside of it; b) the factorisation of the matrix~\eqref{AKM2} with $q\in {\mathcal W}({\mathbb T})$ and $q^2$ being a rational function free of zeros and poles on ${\mathbb T}$.

\subsubsection{Linear-fractional problem approach}
\label{l-f_p}

In \cite{Kiy13} a constructive method of factorisation was proposed, based on the relation of the homogeneous Riemann--Hilbert boundary value problem (\ref{R0}) and so called linear-fractional problem.

Consider a homogeneous  Riemann--Hilbert boundary value problem
\begin{equation}
\label{eq:l-f_p}
W^{+}(t) = G(t) W^{-}(t), \quad t\in {\mathcal L}, \; 0\in {\mathrm{int}}{\mathcal L},\; \infty\in {\mathrm{ext}} {\mathcal L},
\end{equation}
with a non-singular H\"older continuous $2\times 2$ matrix coefficient
\begin{equation}
G(t) = \left(\begin{array}{cc} g_{11}(t) & g_{12}(t) \\ g_{21}(t) & g_{22}(t) \end{array}\right).
\end{equation}
Denote by $W^{+}(z) = \left(w_1^{+}(z),w_2^{+}(z)\right)^T$,  $W^{-}(z) = \left(w_1^{-}(z),w_2^{-}(z)\right)^T$ the components of the solution to \eqref{eq:l-f_p}
and introduce the following functions
\begin{equation}
\label{l-f_p2}
\Phi^{+}(z) = \frac{w_2^{+}(z)}{w_1^{+}(z)},\quad \Phi^{-}(z) = \frac{w_2^{-}(z)}{w_1^{-}(z)}.
\end{equation}
These functions are the components of a piece-wise meromorphic solution to a so-called linear-fractional boundary value problem
\begin{equation}
\label{l-f_p3}
g_{11}(t) \Phi^{+}(t) - g_{22}(t) \Phi^{-}(t) + g_{12}(t) \Phi^{+}(t) \Phi^{-}(t) = g_{21}(t).
\end{equation}
The approach in \cite{Kiy13} (see also \cite{Kiy12}) is based on different possibilities arising from the study of solvability
and representation of the solution to (\ref{l-f_p3}). These
situations are listed in \cite{Kiy12,Kiy13}. The representation of the
solution to the corresponding factorisation problem is given for each
special case. In \cite{Kiy15} this approach is applied to the study of
the Riemann--Hilbert boundary value problem with a $3\times 3$ matrix coefficient $G(t)$.

\subsection{Rawlins--Williams class of matrix Wiener--Hopf problems}
\label{sec:other}

In acoustics and electromagnetism the matrix kernel \(A
(\alpha)\) sometimes is a function only of $\gamma(\alpha)=(k^2- \alpha ^2)^{1/2}$. This motivates the class considered by Rawlins and Williams in~\cite{Raw}
 \[A (\alpha) =\left(
 \begin{array}{cc}
  F(\gamma) &  G (\gamma)F(\gamma)  \\
 H(\gamma) & -G(\gamma)H(\gamma)
 \end{array} \right),\]
where \( F,G\) and $H$ are analytic functions (except possibly when
$\gamma=0$).

The way that this class of matrices is solved is by first transforming them to a matrix Riemann--Hilbert problem on a half line (along a cut emanating from one of the branch-points of $\gamma(\alpha)$). It then transpires that this system can be decoupled into two scalar Riemann--Hilbert problems by Hurd's method~\cite{Hur76}, which can be solved explicitly.

Another class of matrices was considered by D. S. Jones in~\cite{Jones84}:
\[\mathbf{C}  =\left(
 \begin{array}{cc}
  f &  ge_1  \\
 ge_2 & f- 2ge_3
 \end{array} \right),\]
 where  $f(s)$, $g(s)$ are analytic in the strip and $e_1$, $e_2$ and
 $e_3$ are analytic functions in the whole complex plane.
 A commutative factorisation is constructed. There are also extensions
 to non-commutative factorisation by pre-multiplying by  analytic
 matrices. Although it is not obvious, the class of matrices considered by
Rawlins and Williams can, in fact, be reduced to Jones' class~\cite{Jones84}.

For a review of other classes which originate from
applications see~\cite{daniele2014wiener}.

\section{A list of approximate procedures}
\label{sec:approx}

Due to the lack of a general exact constructive factorisation procedure, there has been considerable interest in approximate methods for
Wiener--Hopf matrix factorisation. Numerous approaches have been proposed
in the literature \cite{assym_2, Dan+Lom07,My1, Nigel_cyl,
  Abraha_all_pade, gratting, Lorna_Peake, Lorna, Ant_crack,
  Mishuris-Rogosin, owl, Crighton_matching, MyD-K, Mishuris09, Abr_exp, Shanin_13}. We describe some of the more widely applicable classes here.
Most approximate constructive methods rely on approximating a Wiener--Hopf
problem by a class where exact
constructive methods exist. There are also methods which reduce the
Wiener--Hopf problem to a different equation such as Fredholm integral
equation~\cite{Dan+Lom07, Daniele_Lombardi_11, Daniele_Lombardi_16} but this is outside the scope of this survey.

Suppose there are two functions
\begin{equation}
  \label{eq:perturb}
  K(\alpha)=K_{-}(\alpha)K_{+}(\alpha), \quad \text{and} \quad
\bar{K}(\alpha)=\bar{K}_{-}(\alpha)\bar{K}_{+}(\alpha).
\end{equation}
\begin{Question}
 Is it true
that if \(\bar{K}(\alpha)\) approximates \(K(\alpha)\) than
\(\bar{K}_{\pm}(\alpha)\) approximates \(K_{\pm}(\alpha)\)?
\end{Question}
As stated
the question is too vague to have a definitive answer, but nevertheless this kind
of question motivates the search for approximate solutions.
For scalar functions and \(L_p\) norm  if \(|K(\alpha)-\bar{K}(\alpha)|_{p} \le \epsilon_p\)
there are bounds for \(|K_{\pm}(\alpha)-\bar{K}_{\pm}(\alpha)|_{p}\) in
terms of the \(\epsilon_p\) and computable quantities of
\(K(\alpha)\)~\cite{My1}. The general answer for a matrix valued function
\(K(\alpha)\) is negative. This can be the case due to a choice of
norm or it can be due to an instability.

The simplest example of instability is obtained by mapping an example
 \cite{Spit}
from the unit circle  to the real line. Consider a diagonal matrix
function with partial indices [\(1\), \(-1\)].
\begin{equation}
\label{eq:unstab}
\left(
 \begin{array}{cc}

\frac{t-i}{t+i} &  0 \\
0  & \frac{t+i}{t-i}
 \end{array} \right)=\mathbf{I}
 \left(
 \begin{array}{cc}

\frac{t-i}{t+i} &  0 \\
0  & \frac{t+i}{t-i}
 \end{array} \right)\mathbf{I}.
\end{equation}
Perturbing the matrix we have:
\begin{equation}
\label{eq:unstab_1}
\left(
 \begin{array}{cc}

\frac{t-i}{t+i} &  0 \\
\epsilon  & \frac{t+i}{t-i}
 \end{array} \right)=
\left(
 \begin{array}{cc}

1 &  \frac{t-i}{t+i} \\
0  & \epsilon
 \end{array} \right)
 \mathbf{I}
\left(
 \begin{array}{cc}

0 &  -1/\epsilon \\
1  & \frac{t+i}{\epsilon(t-i)}
 \end{array} \right).
\end{equation}
This example demonstrates that a small perturbation can change the
factors by an arbitrary amount and
can also change the partial indices (from \(\{1\), \(-1\}\) to \(\{0\),
\(0\}\)). This is significant because the partial indices are uniquely defined.
Note that the sum of the partial indices remains the same; this
is true in general, by considering the determinant of both sides
it is reduced to scalar factorisation. For scalar factorisation of
function \(f\)
the index is the winding number of the curve \((\text{Re } f(t),
\text{Im } f(t)) \) \(t \in \mathbb{R}\)  hence is stable under perturbations.  The partial indices are linked to the growth at
infinity of the Wiener--Hopf factorisation, see~\cite{Mer}.

The following surprising theorem provides the necessary and sufficient
conditions for the partial indecies to remain the same under small
enough perturbations.

\begin{thm}[Gohberg -- Krein]\cite{Spit_book}
The system \(\kappa_1 \ge \dots \ge \kappa_n\) of partial indices is stable if and only if:
\[\kappa_1-\kappa_n\le 1.\]
\end{thm}
In fact, this condition is enough to ensure the stability of factors
in the Wiener norm.

\begin{thm}[Shubin]\cite{Spit_book}
Assume the matrix function \(\mathbf{G}\) to have a Wiener--Hopf factorisation and the tuple of its partial indices to be stable. Then, for every \(\epsilon>0\)
there exists a \(\delta>0\) such that, for \(||\mathbf{F}-\mathbf{G}||<\delta\), the matrix function \(\mathbf{F}\) admits a factorisation in which \(||\mathbf{F}_{\pm}-\mathbf{G}_{\pm}||<\epsilon\).
\end{thm}
An obstacle in using this result in applications is that one cannot in
general
determine the partial indices without constructing the factorisation.
Thus, in general it is impossible to prove that some approximate
factorisation of a matrix function is a good approximation. There are
some specific results for Khrapkov--Daniele matrix functions~\cite{MyD-K} and for small perturbations to
the identity matrix~\cite{Mishuris-Rogosin}. This will be discussed in
more detail in the sections below. Interestingly, all three
methods below reduce the matrix factorisation problem to a series
of scalar additive Wiener--Hopf splittings. This, in some sense, avoids
problems with instabilities of matrix factorisation with unstable
indices.

\subsection{Rational approximation}
\label{sec:Rat}

It is attractive to consider a rational matrix \(\bar{K}(\alpha)\) in \eqref{eq:perturb} since
its factorisation can be computed exactly (Section~\ref{constructive}\ref{rational}), without the need for a Cauchy type integral~\eqref{eq:R6} representation. Approximate rational solutions were considered early on~\cite{Koiter} but
they were mainly constructed using \emph{ad hoc} observations
\cite{bookWH}*{Ch. 4.5}. In 2000 a systematic way of
approximating the Wiener--Hopf equations was developed using a two point Pad\'{e} approximation
(the two
points being 0 and \(\infty\) on the real axis)~\cite{Pade}. The power
of the method came from the fact that not the whole matrix was
approximated but only specific parts and  then pole removal method was
applied (Section~\ref{rational}). Since then it proved popular and
found applications in
different branches of mathematics~\cite{Nigel_poles,Abraha_all_pade, Faranosov17}  including finance
\cite{finance}.

From the last paragraph it is clear that it is desirable to
approximate a matrix function, the whole or various parts, by rational
functions \(\bar{K}(\alpha)\).
A natural question to ask: which \(K(\alpha)\) can be
approximated by rational matrix functions \(\bar{K}(\alpha)\)? For this question  we have to
specify in what sense are \(K(\alpha)\) and \(\bar{K}(\alpha)\)
close. It is not an easy task to decide what the correct norm is to
consider, but good candidates are \(L_p\) norms or Sobolev
spaces~\cite{Speck_19}. We will review some known facts about approximation by
rational functions, mostly by Pad\'{e} approximants.
Functions with branch cuts (multi-valued) cannot be ``fully'' approximated by rational
functions which are single-valued. But for a function analytic at
infinity the next best result is true, the maximal domain of \(K(\alpha)\) where the function has a single-valued branch is the domain of convergence of the diagonal Pad\'{e} approximants
for \(K(\alpha)\). In fact, most of the
poles tend to the boundary of the domain of convergence
and lie on the system of cuts that makes the function
single-valued~\cite{Stahl97, Aptekarev2015}. Numerical experimentation
has confirmed this holds true for even small degree of rational
approximation and, what is more, often the interlacing of poles
and zeros occurs on the branch
cut~\cite{Pade_numeric_test,My1}. This simulates the behaviour of the
branch cut; the values on either side of the cut made by poles and zeros
has a jump that tends to the exact value as the rational approximation degree is increased~\cite{Pade_Lushnikov}.

The practical implementation of the rational approximations has now
been well developed. In particular there are numerous algorithms on
\emph{Chebfun} (MATLAB). One of the latest additions in \emph{Chebfun} is the very fast AAA rational approximation~\cite{AAA_Trefethen}.
It is also useful to construct a  mapping of the real line to the unit
interval since most existing algorithms are for bounded intervals~\cite{My1}.
The main difficulty with rational approximations arise when a branch cut
goes though infinity as in \(\gamma (\alpha)=\sqrt{\alpha^2-k^2}\) (which is frequently
encountered in diffraction problems). The problem arises because the
real line crosses the branch cut at infinity. Hence  rational approximation can no
longer be accurate on the whole contour (the real line). It can still be very accurate
on a bounded domain and hence rational approximations of \(\gamma
(\alpha)\) are used in acoustics far-field calculations~\cite{Lorna}.



\subsection{Asymptotic Wiener--Hopf factorisation}
\label{sec:Asy}

Let $M$ be a matrix function analytic around the real axis $\mathbb{R}$.
Let $I$ be the identity matrix and $\varepsilon \in
[0,1]$. We can write
\begin{equation}
	M(x) =I + \varepsilon\, G(x),
\end{equation}
and then the asymptotic algorithm produces a good approximate  Wiener--Hopf
factorisation if \(\varepsilon\, G(x)\) is small compared to
\(I\)~\cite{Mishuris-Rogosin,Crighton_matching}.

Asymptotic Wiener--Hopf factorisation looks for an approximate factorisation of the form
\begin{equation}
I + \varepsilon\, G \approx \bigg(I + \sum_{i \in \mathbb{N^*}}\varepsilon^i\, G^-_{i-1} \bigg)\bigg(I + \sum_{i \in \mathbb{N}^*}\varepsilon^i\, G^+_{i-1} \bigg).
 \end{equation}
 We define the error $\Delta_j$ at a step $j$ as the term neglected in the approximation:
\begin{equation}
   \Delta_j = I + \varepsilon\, G - S_j^-S_j^+ \, ,
\; \text{where} \;
S_j^- = I + \sum_{i=1}^{j}\varepsilon^i\, G^-_{i-1},
\; \text{and} \;
S_j^+ =  I + \sum_{i=1}^{j}\varepsilon^i\, G^+_{i-1}.
\end{equation}

For the first step, we are looking for $G^-_{0}$ and $G^+_{0}$, respectively analytic on the lower half-plane and the upper half-plane, so that
\begin{equation*}
\begin{split}
   I + \varepsilon\, G &\approx (I + \varepsilon\, G^-_{0} )(I + \varepsilon\, G^+_{0} ),\\
	&= I + \varepsilon\,(G^-_{0} + G^+_{0}) + \varepsilon^2\,(G^-_{0}G^+_{0}).
\end{split}
\end{equation*}
Hence, at the first order of $\varepsilon$, we have the Wiener--Hopf additive splitting:
\begin{equation}\label{eq:G0}
G = G^-_{0} + G^+_{0} \, .
\end{equation}
The error is
\begin{equation}
   \Delta_1 = (I + \varepsilon\, G) - (I + \varepsilon\, G^-_{0} )(I + \varepsilon\, G^+_{0} )\\
	= \varepsilon^2\,G^-_{0}G^+_{0}.
\end{equation}
For the second order of $\varepsilon$:
\begin{equation*}
\begin{split}
   I + \varepsilon\, G &\approx (I + \varepsilon\, G^-_{0} + \varepsilon^2\, G^-_{1})(I + \varepsilon\, G^+_{0} + \varepsilon^2\, G^+_{1}),\\
	&= I + \varepsilon\,(G^-_{0} + G^+_{0}) + \varepsilon^2\,(G^-_{1}+G^+_{1}+G^-_{0}G^+_{0})+ \varepsilon^3\,(G^-_{1}G^+_{0}+G^-_{0}G^+_{1})+\varepsilon^4\,G^-_{1}G^+_{1},
\end{split}
\end{equation*}
where the term in $\varepsilon^2\ $ should be equal to zero.
Therefore, we can deduce $G^-_{1}$, $G^+_{1}$ and thus $\Delta_2$:
\begin{equation}\label{eq:step2}
 G_{1}=G^-_{1}+G^+_{1} = -G^-_{0}G^+_{0} \,, \qquad
  \Delta_2 = \varepsilon^3\,(G^-_{1}G^+_{0}+G^-_{0}G^+_{1})+\varepsilon^4\,G^-_{1}G^+_{1} \, .
\end{equation}
In the
same fashion arbitrary order approximations are constructed and, under some conditions, its convergence can be established \cite{Mishuris-Rogosin} and used in \cite{Maurya19} for the factorisation on the unit circle. Furthermore, this technique can be extended for the case of a set of stable partial partial indices \cite{Mishuris09} or, in some cases, even for unstable sets \cite{MR_2018}.

\subsection{Neglecting coupling, and iteration}
\label{sec:Neg}

The main difficulty in solving matrix Wiener--Hopf equation arises
from the fact that there are {\it coupled} scalar equations. In some
situations the coupling is small and a reasonable approximation can be
obtained by neglecting it. In other words the original matrix is close
to being diagonal or triangular. More often the coupling is too strong
to give a desired approximation see Section 4.4 in~\cite{bookWH} or \cite{Kobayashi91}, nevertheless the same ideas can still be utilised.
One way is through the use of iterative methods, which in the context of
Wiener--Hopf equations were proposed early on~\cite{old_iter,BELINSKII_iter}. One of the reasons that they did not gain popularity in the past is that they are  computationally intensive. But recent advances in computing Cauchy-type transforms, and significant computational resource, have now made it a practical approach~\cite{olver_computing_2011,olver_change,trogdon_riemannhilbert_2016}. In
particular it has been applied to a class of triangular matrix functions containing exponential
factors~\cite{My_iter}. The aim is to find functions \(\Phi_-^{(0)}(\alpha)\),
\(\Phi_-^{(L)} (\alpha)\),
\(\Psi_+^{ (0)}(\alpha)\)
and \(\Psi_+^{ (L)}(\alpha)\)
analytic in respective half-planes, satisfying the following
relationship
\begin{equation}
  \label{eq:main}
\begin{pmatrix}
 \Phi_-^{(0)}(\alpha) \\
 \Phi_-^{(L)} (\alpha)
 \end{pmatrix}=\begin{pmatrix}
  A(\alpha)&  B(\alpha)e^{i \alpha L} \\
 C(\alpha)e^{-i \alpha L} & 0
 \end{pmatrix}\begin{pmatrix}
 \Psi_+^{ (0)}(\alpha) \\
 \Psi_+^{ (L)}(\alpha)
 \end{pmatrix}+ \begin{pmatrix}
 f_1(\alpha)\\
f_2 (\alpha)
 \end{pmatrix},
\end{equation}
on the strip \(\mathcal H \). The functions \(A(\alpha)\),
\(B(\alpha)\) and \(C(\alpha)\) are known and \(L\) is a positive
constant. The method has found applications in
aeroacoustics~\cite{my+Lorna}, crack propagation
problems~\cite{Mish_crack_iter}.
It has been extended to $n$
dimensional matrix problems and applied to scattering from multiple
co-linear plates~\cite{iter_n}. We would also like to note that there is a similarity of the
above iteration methods with a concept in diffraction called the
Schwarzschild series~\cite{Schwarzschild1901,ufimtsev2003theory, Shanin_S_3,Sautbekov_17}.
This illustrates a useful viewpoint: many  problems which give rise to a matrix Wiener--Hopf equation can be
thought of as an interaction of a number of distinct problems, each of which would lead to a scalar Wiener--Hopf equation~\cite{Abrahams_82}.

\section{Related methods and open problems}

To solve Wiener--Hopf type problems which
arise in applications, the following are desirable
\begin{enumerate}
\item A systematic way of determining if the posed problem is of
  Wiener--Hopf type and has a unique solution.
\item A justified and executable algorithm capable of verifying whether a given matrix possesses a canonical factorisation (all partial indices are equal to zero) or not.
\item A routine method of deriving and transforming the original problem into a Wiener--Hopf equation with regularised matrix kernel (allowing canonical factorisation).
  \item Criteria for determining which class of matrix Wiener--Hopf
    equation the given equation belongs to.
    \item A unified constructive method for an exact or an approximate factorisation.
    \end{enumerate}
 For problems of scalar Wiener--Hopf type, the above points can be considered
 relatively complete.
 In the case of a classical matrix Wiener--Hopf equations the first one is
 straightforward while the remaining points continue to offer future avenues of
 research. It is possible to develop alternative methods for solving
 equations of Wiener--Hopf type, for example based on Fredholm integral equation theory~\cite{Dan+Lom07, Book_Daniele20_1, Book_Daniele20_2}, and also more theoretical operator based approaches to convolution integral equations \cite{Sakhnovich_15, Speck_19}.

 There  are equations of Wiener--Hopf type which have not been
 reviewed here.
 For example, PDEs on wedge domains are not immediately in a standard class, but can be stated
as generalised Wiener--Hopf equations if appropriate transformations/mapping
are considered~\cite{Nethercote2020a, BELINSKII_elasticright}. The difficulty with wedge domains is that the Fourier
transform is not the natural operator to use, hence a mapping is
required. Alternatively, the Mellin transform~\cite{Sneddon95} could
be employed if the respective differential operators are of the same
homogeneity, and a combination of the Mellin and Fourier transform could be used for multilayered-multiwedged domains~\cite{MishurisOlesiak95, Mishris96, Mishuris97}. There is also a link to the unified transform method which provides
a systematic way of deriving a natural transform for some polygonal
domains~\cite{FokSpen12}. A natural but difficult generalisation is to
consider~\eqref{eq:WH1} in more variables (a double integral)~\cite{Assier20,Assier21} leading to interesting
Wiener--Hopf related methods in multi-variate complex
analysis~\cite{Assier2018, Assier_Shanin}. A different direction is when the kernel \eqref{eq:WH1} has support on the half-line but is not of convolution
type (the equation is referred to as being of Hammerstein type)~\cite{non-con}. There is also a wealth of literature on spectral properties of
Toeplitz matrices/operators \cite{Bottcher_spectrum_book05,
  Grudsky_08, Farenick_96, Bottcher_11} which is a
related topic but outside the scope of this review.

Interest in Wiener--Hopf techniques has been steady
ever since the 1970s, judging by the number of papers containing ``Wiener--Hopf''
on {\it mathscinet}, with a total around 3,000 to date (it is likely to
be a significant underestimate since many of the references for this
review either do not contain ``Wiener--Hopf'' in the title, or are published in engineering and physics journals). The related
Riemann--Hilbert boundary value method has also been popular with around 4,600
publications; with a small but steady increase in number of
publication each year overall.  The Wiener--Hopf method could be viewed through many lenses, has many related problems, arises in numerous applications and still has many interesting open problems.
It can be expected that further work in coming years will resolve some of the
long-standing questions related to, and limitations of, the Wiener--Hopf method, and hence extend its relevance and applicability still further.






\subsection{Funding}

The authors would like to thank the Isaac
Newton Institute for Mathematical Sciences, Cambridge UK, for support and hospitality during
the programme ``Bringing pure and applied analysis together via the Wiener--Hopf technique,
its generalisations and applications'' where some of the work on this article was undertaken
(supported by EPSRC grant no EP/R014604/1).  A.V.K. is supported by  Royal Society
Dorothy Hodgkin Research Fellowship and Dame Kathleen Ollerenshaw
Fellowship. S.V.R. is partially supported by the Belarusian Fund for
Fundamental Research through grant F20MS-083. G.M. is supported by the Royal Society Wolfson Research Merit Award and Ser Cymru Future Generations Industrial Fellowship.

\subsection{Acknowledgement}
  A.V.K. would like to thank Yuri Antipov, Andrey Shanin and
  Boris Belinskiy for interesting and useful discussions. G.M and S.V.R. would like to thank Ilya Spitkovsky for relevant discussions.


\bibliographystyle{RS}
\bibliography{review}

\def\cprime{$'$}
\begin{bibdiv}
\begin{biblist}

\bib{NSR_20}{article}{
      author={Abrahams, David},
      author={Huang, Xun},
      author={Kisil, Anastasia},
      author={Mishuris, Gennady},
      author={Nieves, Michael},
      author={Rogosin, Sergei},
      author={Spitkovsky, Ilya},
       title={{Reinvigorating the Wiener-Hopf technique in the pursuit of
  understanding processes and materials}},
        date={202009},
        ISSN={2095-5138},
     journal={National Science Review},
      volume={8},
      number={2},
  eprint={https://academic.oup.com/nsr/article-pdf/8/2/nwaa225/38626007/nwaa225.pdf},
         url={https://doi.org/10.1093/nsr/nwaa225},
        note={nwaa225},
}

\bib{Abrahams_82}{article}{
      author={Abrahams, I.~D.},
       title={Scattering of sound by large finite geometries},
        date={1982},
        ISSN={0272-4960},
     journal={IMA J. Appl. Math.},
      volume={29},
      number={1},
       pages={79\ndash 97},
         url={https://doi-org.manchester.idm.oclc.org/10.1093/imamat/29.1.79},
      review={\MR{667783}},
}

\bib{Elastic_half}{article}{
      author={Abrahams, I.~D.},
       title={Radiation and scattering of waves on an elastic half-space; a
  non-commutative matrix {W}iener-{H}opf problem},
        date={1996},
        ISSN={0022-5096},
     journal={J. Mech. Phys. Solids},
      volume={44},
      number={12},
       pages={2125\ndash 2154},
  url={http://0-dx.doi.org.pugwash.lib.warwick.ac.uk/10.1016/S0022-5096(96)00064-6},
      review={\MR{1423542 (97i:73044)}},
}

\bib{Abr_exp}{article}{
      author={Abrahams, I.~D.},
      author={Wickham, G.~R.},
       title={General {W}iener-{H}opf factorization of matrix kernels with
  exponential phase factors},
        date={1990},
        ISSN={00361399},
     journal={SIAM Journal on Applied Mathematics},
      volume={50},
      number={3},
       pages={819\ndash 838},
         url={http://www.jstor.org/stable/2101888},
}

\bib{Pade}{article}{
      author={Abrahams, I.~David},
       title={The application of {P}ad\'e approximations to {W}iener-{H}opf
  factorization},
        date={2000},
        ISSN={0272-4960},
     journal={IMA J. Appl. Math.},
      volume={65},
      number={3},
       pages={257\ndash 281},
         url={http://dx.doi.org/10.1093/imamat/65.3.257},
      review={\MR{1806416 (2001k:45010)}},
}

\bib{Abrahams_02}{incollection}{
      author={Abrahams, I.~David},
       title={On the application of the {W}iener-{H}opf technique to problems
  in dynamic elasticity},
        date={2002},
      volume={36},
       pages={311\ndash 333},
  url={https://doi-org.manchester.idm.oclc.org/10.1016/S0165-2125(02)00027-6},
        note={Dedicated to Jan D. Achenbach on the occasion of his 65th
  birthday},
      review={\MR{1950990}},
}

\bib{Abraha_all_pade}{article}{
      author={Abrahams, I.~David},
      author={Davis, Anthony M.~J.},
      author={Llewellyn~Smith, Stefan~G.},
       title={Matrix {W}iener-{H}opf approximation for a partially clamped
  plate},
        date={2008},
        ISSN={0033-5614},
     journal={Quart. J. Mech. Appl. Math.},
      volume={61},
      number={2},
       pages={241\ndash 265},
         url={http://dx.doi.org/10.1093/qjmam/hbn004},
      review={\MR{2414435 (2009d:74047)}},
}

\bib{Adu93}{article}{
      author={Adukov, V.~M.},
       title={{Wiener-Hopf factorization of meromorphic matrix-functions}},
        date={1993},
     journal={St. Petersburg Math. J.},
      volume={4},
      number={1},
       pages={51\ndash 69},
}

\bib{Adu99}{article}{
      author={Adukov, V.~M.},
       title={{Factorization of analytic matrix-valued functions}},
        date={1999},
     journal={Theoret. and Math. Phys.},
      volume={118},
      number={3},
       pages={255\ndash 263},
}

\bib{Adu13}{article}{
      author={Adukov, V.~M.},
       title={{On Wiener-Hopf factorization of functionally-commutative
  matrix-functions}},
        date={2013},
     journal={Vestnik South Ural State University, Ser. Math. Mech. Phys.
  (Russian)},
      volume={5},
      number={2},
       pages={6\ndash 12},
}

\bib{Adu18}{article}{
      author={Adukov, V.~M.},
       title={{Algorithm of polynomial factorization and its implementation in
  Maple}},
        date={2018},
     journal={Vestnik South Ural State University, Ser. Math. Mech. Phys.
  (Russian)},
      volume={11},
      number={4},
       pages={110\ndash 122},
}

\bib{Adukov_20}{article}{
      author={Adukov, Victor~M.},
      author={Mishuris, Gennady},
      author={Rogosin, Sergei~V.},
       title={Exact conditions for preservation of the partial indices of a
  perturbed triangular {$2\times 2$} matrix function},
        date={2020},
        ISSN={1364-5021},
     journal={Proc. A.},
      volume={476},
      number={2237},
       pages={20200099, 19},
         url={https://doi-org.manchester.idm.oclc.org/10.1098/rspa.2020.0099},
      review={\MR{4111977}},
}

\bib{Akto92}{article}{
      author={Aktosun, T.},
      author={Klaus, M.},
      author={van~der Mee, C.},
       title={{Explicit Wiener-Hopf factorization for certain non-rational
  matrix functions}},
        date={1992},
     journal={Integral Equations Operator Theory},
      volume={15},
       pages={879\ndash 900},
}

\bib{Albani_meta}{article}{
      author={Albani, Matteo},
      author={Capolino, Filippo},
       title={Wave dynamics by a plane wave on a half-space metamaterial made
  of plasmonic nanospheres: a discrete wiener--hopf formulation},
        date={2011Sep},
     journal={J. Opt. Soc. Am. B},
      volume={28},
      number={9},
       pages={2174\ndash 2185},
         url={http://josab.osa.org/abstract.cfm?URI=josab-28-9-2174},
}

\bib{Antipov99}{article}{
      author={Antipov, Y.~A.},
       title={An exact solution of the {$3$}-{D}-problem of an interface
  semi-infinite plane crack},
        date={1999},
        ISSN={0022-5096},
     journal={J. Mech. Phys. Solids},
      volume={47},
      number={5},
       pages={1051\ndash 1093},
  url={https://doi-org.manchester.idm.oclc.org/10.1016/S0022-5096(98)00096-9},
      review={\MR{1687091}},
}

\bib{Antipov_09}{article}{
      author={Antipov, Y.~A.},
       title={Subsonic semi-infinite crack with a finite friction zone in a
  bimaterial},
        date={2009},
        ISSN={0022-5096},
     journal={J. Mech. Phys. Solids},
      volume={57},
      number={12},
       pages={1934\ndash 1957},
  url={https://doi-org.manchester.idm.oclc.org/10.1016/j.jmps.2009.08.006},
      review={\MR{2583608}},
}

\bib{Antipov_18}{article}{
      author={Antipov, Y.~A.},
       title={Subsonic frictional cavitating penetration of a thin rigid body
  into an elastic medium},
        date={2018},
        ISSN={0033-5614},
     journal={Quart. J. Mech. Appl. Math.},
      volume={71},
      number={2},
       pages={221\ndash 243},
         url={https://doi-org.manchester.idm.oclc.org/10.1093/qjmam/hby003},
      review={\MR{3800008}},
}

\bib{Antipov03}{article}{
      author={Antipov, Y.~A.},
      author={Chuang, T.-J.},
      author={Gao, H.},
       title={On the integro-differential equation associated with diffusive
  crack growth theory},
        date={2003},
        ISSN={0033-5614},
     journal={Quart. J. Mech. Appl. Math.},
      volume={56},
      number={2},
       pages={289\ndash 310},
         url={https://doi-org.manchester.idm.oclc.org/10.1093/qjmam/56.2.289},
      review={\MR{1982979}},
}

\bib{Antipov_DK}{article}{
      author={Antipov, Y.~A.},
      author={Mkhitaryan, S.~M.},
       title={A crack induced by a thin rigid inclusion partly debonded from
  the matrix},
        date={2017},
        ISSN={0033-5614},
     journal={Quart. J. Mech. Appl. Math.},
      volume={70},
      number={2},
       pages={153\ndash 185},
         url={https://doi-org.manchester.idm.oclc.org/10.1093/qjmam/hbx003},
      review={\MR{3655936}},
}

\bib{AntSil02}{article}{
      author={Antipov, Y.~A.},
      author={Silvestrov, V.~V.},
       title={{Factorization on a Riemann surface in scattering theory}},
        date={2002},
     journal={Quart. J. Mech. Appl. Math.},
      volume={55},
       pages={607\ndash 654},
}

\bib{AntSil04}{article}{
      author={Antipov, Y.~A.},
      author={Silvestrov, V.~V.},
       title={Second-order functional-difference equations. {II}. {S}cattering
  from a right-angled conductive wedge for {$E$}-polarization},
        date={2004},
        ISSN={0033-5614},
     journal={Quart. J. Mech. Appl. Math.},
      volume={57},
      number={2},
       pages={267\ndash 313},
         url={https://doi-org.manchester.idm.oclc.org/10.1093/qjmam/57.2.267},
      review={\MR{2068407}},
}

\bib{Ant_crack}{article}{
      author={Antipov, Y.~A.},
      author={Smirnov, A.~V.},
       title={Fundamental solution and the weight functions of the transient
  problem on a semi-infinite crack propagating in a half-plane},
        date={2016},
        ISSN={1521-4001},
     journal={ZAMM - Journal of Applied Mathematics and Mechanics},
      volume={96},
      number={10},
       pages={1156\ndash 1174},
         url={http://dx.doi.org/10.1002/zamm.201400272},
}

\bib{Ant14}{article}{
      author={Antipov, Y.A.},
       title={{The Baker-Akhiezer function and factorization of the
  Chebotarev-Khrapkov matrix}},
        date={2014},
     journal={Letters in Mathematical Physics},
      volume={104},
       pages={1365\ndash 1384},
}

\bib{Aptekarev2015}{article}{
      author={Aptekarev, Alexander},
      author={Yattselev, Maxim},
       title={Pad approximants for functions with branch points strong
  asymptotics of nuttall stahl polynomials},
    language={eng},
        date={2015},
        ISSN={0001-5962},
     journal={Acta Mathematica},
      volume={215},
      number={2},
       pages={217\ndash 280},
}

\bib{Assier2018}{article}{
      author={Assier, R.~C.},
      author={Shanin, A.~V.},
       title={Diffraction by a quarter{\textendash}plane. analytical
  continuation of spectral functions},
        date={2018-12},
     journal={QJMAM},
      volume={72},
      number={1},
       pages={51\ndash 86},
         url={https://doi.org/10.1093/qjmam/hby021},
}

\bib{Assier20}{article}{
      author={Assier, Rapha\"{e}l~C.},
      author={Abrahams, I.~David},
       title={On the asymptotic properties of a canonical diffraction
  integral},
        date={2020},
        ISSN={1364-5021},
     journal={Proc. A.},
      volume={476},
      number={2242},
       pages={150\ndash 165},
         url={https://doi-org.manchester.idm.oclc.org/10.1098/rspa.2020.0150},
      review={\MR{4176312}},
}

\bib{Assier21}{article}{
      author={Assier, Raphael~C.},
      author={Abrahams, I.~David},
       title={A {S}urprising {O}bservation in the {Q}uarter-{P}lane
  {D}iffraction {P}roblem},
        date={2021},
        ISSN={0036-1399},
     journal={SIAM J. Appl. Math.},
      volume={81},
      number={1},
       pages={60\ndash 90},
         url={https://doi-org.manchester.idm.oclc.org/10.1137/19M1258785},
      review={\MR{4199740}},
}

\bib{Assier_Shanin}{article}{
      author={Assier, Raphael~C.},
      author={Shanin, Andrey~V.},
       title={Analytical continuation of two-dimensional wave fields},
        date={2021},
     journal={Proceedings of the Royal Society A: Mathematical, Physical and
  Engineering Sciences},
      volume={477},
      number={2245},
       pages={20200681},
  eprint={https://royalsocietypublishing.org/doi/pdf/10.1098/rspa.2020.0681},
  url={https://royalsocietypublishing.org/doi/abs/10.1098/rspa.2020.0681},
}

\bib{Lorna}{article}{
      author={Ayton, Lorna~J},
       title={Acoustic scattering by a finite rigid plate with a poroelastic
  extension},
        date={2016},
     journal={Journal of Fluid Mechanics},
      volume={791},
       pages={414\ndash 438},
}

\bib{Lorna_Peake}{article}{
      author={Ayton, Lorna~J.},
      author={Peake, N.},
       title={On high-frequency noise scattering by aerofoils in flow},
        date={201311},
        ISSN={1469-7645},
     journal={Journal of Fluid Mechanics},
      volume={734},
       pages={144\ndash 182},
         url={http://journals.cambridge.org/article_S0022112013004771},
}

\bib{Ayzenberg14}{article}{
      author={Ayzenberg-Stepanenko, Mark},
      author={Mishuris, Gennady},
      author={Slepyan, Leonid},
       title={Brittle fracture in a periodic structure with internal potential
  energy. {S}pontaneous crack propagation},
        date={2014},
        ISSN={1364-5021},
     journal={Proc. R. Soc. Lond. Ser. A Math. Phys. Eng. Sci.},
      volume={470},
      number={2167},
       pages={20140121, 20},
         url={https://doi-org.manchester.idm.oclc.org/10.1098/rspa.2014.0121},
      review={\MR{3209223}},
}

\bib{Craster_ice}{article}{
      author={Balmforth, N.~J.},
      author={Craster, R.~V.},
       title={Ocean waves and ice sheets},
        date={1999},
        ISSN={0022-1120},
     journal={J. Fluid Mech.},
      volume={395},
       pages={89\ndash 124},
         url={https://doi.org/10.1017/S0022112099005145},
      review={\MR{1736494}},
}

\bib{Bankov2008}{article}{
      author={Bankov, S.~E.},
       title={Diffraction of electromagnetic waves by the boundary of a
  semi-infinite metamaterial},
        date={2008Jan},
        ISSN={1555-6557},
     journal={Journal of Communications Technology and Electronics},
      volume={53},
      number={1},
       pages={15\ndash 25},
         url={https://doi.org/10.1134/S1064226908010026},
}

\bib{BELINSKII_iter}{article}{
      author={Belinskii, B.P.},
       title={The integral equations of stationary problems on the diffraction
  of short waves at obstacles of the segment type},
        date={1973},
        ISSN={0041-5553},
     journal={USSR Computational Mathematics and Mathematical Physics},
      volume={13},
      number={2},
       pages={125 \ndash  140},
  url={http://www.sciencedirect.com/science/article/pii/0041555373901365},
}

\bib{BELINSKII_elasticright}{article}{
      author={Belinskii, B.P.},
      author={Kouzov, D.P.},
      author={Chel'tsova, V.D.},
       title={On acoustic wave diffraction by plates connected at a right
  angle: Pmm vol. 37, n2, 1973, pp. 291-299},
        date={1973},
        ISSN={0021-8928},
     journal={Journal of Applied Mathematics and Mechanics},
      volume={37},
      number={2},
       pages={273 \ndash  281},
  url={http://www.sciencedirect.com/science/article/pii/0021892873900361},
}

\bib{Bir03}{article}{
      author={Birkhoff, G.~D.},
       title={{The generalized Riemann problem for linear differential
  equations and allied problems for linear difference and $q$-difference
  equations}},
        date={1913},
     journal={Proc. Amer. Acad. Arts and Sci.},
      volume={49},
      number={9},
       pages={521\ndash 568},
}

\bib{Bol90}{article}{
      author={Bolibruch, A.~A.},
       title={{The Riemann--Hilbert problem}},
        date={1990},
     journal={Russian Math. Surveys},
      volume={45},
      number={2},
       pages={1\ndash 47},
}

\bib{Bottcher_11}{article}{
      author={B\"{o}ttcher, Albrecht},
      author={Grudsky, Sergei},
      author={Iserles, Arieh},
       title={Spectral theory of large {W}iener-{H}opf operators with
  complex-symmetric kernels and rational symbols},
        date={2011},
        ISSN={0305-0041},
     journal={Math. Proc. Cambridge Philos. Soc.},
      volume={151},
      number={1},
       pages={161\ndash 191},
  url={https://doi-org.manchester.idm.oclc.org/10.1017/S0305004111000259},
      review={\MR{2801320}},
}

\bib{Bottcher_spectrum_book05}{book}{
      author={B\"{o}ttcher, Albrecht},
      author={Grudsky, Sergei~M.},
       title={Spectral properties of banded {T}oeplitz matrices},
   publisher={Society for Industrial and Applied Mathematics (SIAM),
  Philadelphia, PA},
        date={2005},
        ISBN={0-89871-599-7},
         url={https://doi-org.manchester.idm.oclc.org/10.1137/1.9780898717853},
      review={\MR{2179973}},
}

\bib{Brun14}{article}{
      author={Brun, Michele},
      author={Giaccu, Gian~Felice},
      author={Movchan, Alexander~B.},
      author={Slepyan, Leonid~I.},
       title={Transition wave in the collapse of the san saba bridge},
        date={2014},
        ISSN={2296-8016},
     journal={Frontiers in Materials},
      volume={1},
       pages={12},
         url={https://www.frontiersin.org/article/10.3389/fmats.2014.00012},
}

\bib{Budaev95}{book}{
      author={Budaev, Bair},
       title={Diffraction by wedges},
      series={Pitman Research Notes in Mathematics Series},
   publisher={Longman Scientific \& Technical, Harlow},
        date={1995},
      volume={322},
        ISBN={0-582-10310-X},
      review={\MR{1413497}},
}

\bib{Albani_19}{article}{
      author={Camacho, Miguel},
      author={Hibbins, Alastair~P.},
      author={Capolino, Filippo},
      author={Albani, Matteo},
       title={Diffraction by a truncated planar array of dipoles: a
  {W}iener-{H}opf approach},
        date={2019},
        ISSN={0165-2125},
     journal={Wave Motion},
      volume={89},
       pages={28\ndash 42},
  url={https://doi-org.manchester.idm.oclc.org/10.1016/j.wavemoti.2019.03.004},
      review={\MR{3922189}},
}

\bib{Camara_sur}{article}{
      author={C\^{a}mara, M.~C.},
      author={dos Santos, A.~F.},
      author={dos Santos, Pedro~F.},
       title={Matrix {R}iemann-{H}ilbert problems and factorization on
  {R}iemann surfaces},
        date={2008},
        ISSN={0022-1236},
     journal={J. Funct. Anal.},
      volume={255},
      number={1},
       pages={228\ndash 254},
  url={https://doi-org.manchester.idm.oclc.org/10.1016/j.jfa.2008.01.008},
      review={\MR{2417816}},
}

\bib{Mer}{article}{
      author={C{\^a}mara, M.~C.},
      author={Lebre, A.~B.},
      author={Speck, F.-O.},
       title={Meromorphic factorization, partial index estimates and
  elastodynamic diffraction problems},
        date={1992},
        ISSN={0025-584X},
     journal={Math. Nachr.},
      volume={157},
       pages={291\ndash 317},
         url={http://dx.doi.org/10.1002/mana.19921570124},
      review={\MR{1233066 (94j:47030)}},
}

\bib{Car22}{article}{
      author={Carleman, T.},
       title={{Sur la r\'{e}solution de certaines \'{e}quations int\'{e}grates}
  (french)},
        date={1922},
     journal={Ark. Mat. Astr. Fys.},
      volume={16},
      number={26},
}

\bib{Toeplitz_fast_07}{article}{
      author={Chandrasekaran, S.},
      author={Gu, M.},
      author={Sun, X.},
      author={Xia, J.},
      author={Zhu, J.},
       title={A superfast algorithm for {T}oeplitz systems of linear
  equations},
        date={2007},
        ISSN={0895-4798},
     journal={SIAM J. Matrix Anal. Appl.},
      volume={29},
      number={4},
       pages={1247\ndash 1266},
         url={https://doi-org.manchester.idm.oclc.org/10.1137/040617200},
      review={\MR{2369294}},
}

\bib{Che56}{article}{
      author={Chebotarev, G.~N.},
       title={{Partial indices of the Riemann boundary value problem with a
  triangular matrix of the second order} (russian)},
        date={1956},
     journal={Uspekhi Mat. Nauk},
      volume={11},
      number={3(69)},
       pages={199\ndash 202},
}

\bib{Che56a}{article}{
      author={Chebotarev, G.~N.},
       title={{To the solution in an explicit form of the Riemann boundary
  value problem for a system of $n$ pairs of functions} (russian)},
        date={1956},
     journal={Uch. zapiski Kazan State University. Mathematika},
      volume={116},
      number={4},
       pages={31\ndash 58},
}

\bib{magneto_18}{article}{
      author={Cohen, Roie},
      author={Goldstein, Moshe},
       title={Hall and dissipative viscosity effects on edge magnetoplasmons},
        date={2018Dec},
     journal={Phys. Rev. B},
      volume={98},
       pages={235103},
         url={https://link.aps.org/doi/10.1103/PhysRevB.98.235103},
}

\bib{Atkinson_94}{article}{
      author={Craster, R.~V.},
      author={Atkinson, C.},
       title={Mixed boundary-value problems in nonhomogeneous elastic
  materials},
        date={1994},
        ISSN={0033-5614},
     journal={Quart. J. Mech. Appl. Math.},
      volume={47},
      number={2},
       pages={183\ndash 206},
         url={https://doi-org.manchester.idm.oclc.org/10.1093/qjmam/47.2.183},
      review={\MR{1277145}},
}

\bib{Crighton_matching}{article}{
      author={Crighton, D.~G.},
       title={Asymptotic factorization of {W}iener-{H}opf kernels},
        date={2001},
        ISSN={0165-2125},
     journal={Wave Motion},
      volume={33},
      number={1},
       pages={51\ndash 65},
  url={http://0-dx.doi.org.pugwash.lib.warwick.ac.uk/10.1016/S0165-2125(00)00063-9},
      review={\MR{1804459 (2001j:78012)}},
}

\bib{Croisille99}{book}{
      author={Croisille, Jean-Pierre},
      author={Lebeau, Gilles},
       title={Diffraction by an immersed elastic wedge},
      series={Lecture Notes in Mathematics},
   publisher={Springer-Verlag, Berlin},
        date={1999},
      volume={1723},
        ISBN={3-540-66810-1},
         url={https://doi-org.manchester.idm.oclc.org/10.1007/BFb0092515},
      review={\MR{1740860}},
}

\bib{Crowdy_con}{article}{
      author={Crowdy, Darren~G.},
       title={Uniform flow past a periodic array of cylinders},
        date={2016},
        ISSN={0997-7546},
     journal={Eur. J. Mech. B Fluids},
      volume={56},
       pages={120\ndash 129},
  url={https://doi-org.manchester.idm.oclc.org/10.1016/j.euromechflu.2015.10.003},
      review={\MR{3472576}},
}

\bib{Dan+Lom07}{article}{
      author={Daniele, V.},
      author={Lombardi, G.},
       title={Fredholm factorization of {W}iener-{H}opf scalar and matrix
  kernels},
        date={2007},
        ISSN={1944-799X},
     journal={Radio Science},
      volume={42},
      number={6},
       pages={n/a\ndash n/a},
         url={http://dx.doi.org/10.1029/2007RS003673},
        note={RS6S01},
}

\bib{Dan78}{article}{
      author={Daniele, V.~G.},
       title={{On the factorization of Wiener-Hopf matrices in problems
  solvable with Hurd's method}},
        date={1978},
     journal={IEEE Trans. Antennas and Propagation},
      volume={26},
       pages={614\ndash 616},
}

\bib{daniele2014wiener}{book}{
      author={Daniele, V.G.},
      author={Zich, R.S.},
       title={The {Wiener--Hopf} method in electromagnetics},
      series={Mario Boella Series on Electromagnetism in Information and
  Communication Series},
   publisher={Inst. Engin. Tech., SciTech Publishing},
        date={2014},
        ISBN={9781613530016},
         url={http://books.google.co.uk/books?id=lftHBAAAQBAJ},
}

\bib{Daniele_Lombardi_11}{article}{
      author={Daniele, Vito},
      author={Lombardi, Guido},
       title={The {W}iener-{H}opf solution of the isotropic penetrable wedge
  problem: diffraction and total field},
        date={2011},
        ISSN={0018-926X},
     journal={IEEE Trans. Antennas and Propagation},
      volume={59},
      number={10},
       pages={3797\ndash 3818},
         url={http://dx.doi.org/10.1109/TAP.2011.2163780},
      review={\MR{2893576}},
}

\bib{Daniele_Lombardi_16}{article}{
      author={Daniele, Vito},
      author={Lombardi, Guido},
       title={Arbitrarily oriented perfectly conducting wedge over a dielectric
  half-space: diffraction and total far field},
        date={2016},
        ISSN={0018-926X},
     journal={IEEE Trans. Antennas and Propagation},
      volume={64},
      number={4},
       pages={1416\ndash 1433},
         url={http://dx.doi.org/10.1109/TAP.2016.2524412},
      review={\MR{3488349}},
}

\bib{Book_Daniele20_1}{book}{
      author={Daniele, Vito},
      author={Lombardi, Guido},
       title={Scattering and diffraction by wedges 1: The wiener--hopf
  solution--theory},
   publisher={London: Wiley-ISTE},
        date={2020},
        ISBN={9781119476733},
}

\bib{Book_Daniele20_2}{book}{
      author={Daniele, Vito},
      author={Lombardi, Guido},
       title={Scattering and diffraction by wedges 2: The wiener-hopf solution
  - advanced applications},
   publisher={London: Wiley-ISTE},
        date={2020},
        ISBN={9781119476733},
}

\bib{Duduchava79}{book}{
      author={Duduchava, Roland},
       title={Integral equations in convolution with discontinuous presymbols,
  singular integral equations with fixed singularities, and their applications
  to some problems of mechanics},
   publisher={BSB B. G. Teubner Verlagsgesellschaft, Leipzig},
        date={1979},
        note={With German, French and Russian summaries, Teubner-Texte zur
  Mathematik. [Teubner Texts on Mathematics]},
      review={\MR{571890}},
}

\bib{Pade_Lushnikov}{article}{
      author={Dyachenko, S.~A.},
      author={Lushnikov, P.~M.},
      author={Korotkevich, A.~O.},
       title={Branch cuts of {S}tokes wave on deep water. {P}art {I}:
  {N}umerical solution and {P}ad\'{e} approximation},
        date={2016},
        ISSN={0022-2526},
     journal={Stud. Appl. Math.},
      volume={137},
      number={4},
       pages={419\ndash 472},
         url={https://doi-org.manchester.idm.oclc.org/10.1111/sapm.12128},
      review={\MR{3570991}},
}

\bib{non-con}{article}{
      author={Eggermont, P. P.~B.},
       title={On noncompact {H}ammerstein integral equations and a nonlinear
  boundary value problem for the heat equation},
        date={1992},
        ISSN={0897-3962},
     journal={J. Integral Equations Appl.},
      volume={4},
      number={1},
       pages={47\ndash 68},
         url={https://doi.org/10.1216/jiea/1181075665},
      review={\MR{1160087}},
}

\bib{Spit_piece-wise}{article}{
      author={Ehrhardt, T.},
      author={Spitkovski{\u\i}, I.~M.},
       title={Factorization of piecewise-constant matrix functions, and systems
  of linear differential equations},
        date={2001},
        ISSN={0234-0852},
     journal={Algebra i Analiz},
      volume={13},
      number={6},
       pages={56\ndash 123},
      review={\MR{1883840 (2004h:30045)}},
}

\bib{Ehr_Spec}{article}{
      author={Ehrhardt, Torsten},
      author={Speck, Frank-Olme},
       title={Transformation techniques towards the factorization of
  non-rational {$2\times 2$} matrix functions},
        date={2002},
        ISSN={0024-3795},
     journal={Linear Algebra Appl.},
      volume={353},
       pages={53\ndash 90},
         url={http://dx.doi.org/10.1016/S0024-3795(02)00288-4},
      review={\MR{1918749 (2003f:15021)}},
}

\bib{gratting}{article}{
      author={Erba{\c{s}}, Bar{\i}{\c{s}}},
      author={Abrahams, I.~David},
       title={Scattering of sound waves by an infinite grating composed of
  rigid plates},
        date={2007},
        ISSN={0165-2125},
     journal={Wave Motion},
      volume={44},
      number={4},
       pages={282\ndash 303},
         url={http://dx.doi.org/10.1016/j.wavemoti.2006.11.001},
      review={\MR{2307116 (2007m:76116)}},
}

\bib{Faranosov17}{article}{
      author={Faranosov, Georgy~A.},
      author={Bychkov, Oleg~P.},
       title={Two-dimensional model of the interaction of a plane acoustic wave
  with nozzle edge and wing trailing edge},
        date={2017},
     journal={The Journal of the Acoustical Society of America},
      volume={141},
      number={1},
       pages={289\ndash 299},
      eprint={https://doi.org/10.1121/1.4973855},
         url={https://doi.org/10.1121/1.4973855},
}

\bib{Farenick_96}{article}{
      author={Farenick, Douglas~R.},
      author={Lee, Woo~Young},
       title={Hyponormality and spectra of {T}oeplitz operators},
        date={1996},
        ISSN={0002-9947},
     journal={Trans. Amer. Math. Soc.},
      volume={348},
      number={10},
       pages={4153\ndash 4174},
  url={https://doi-org.manchester.idm.oclc.org/10.1090/S0002-9947-96-01683-2},
      review={\MR{1363943}},
}

\bib{FGK94}{article}{
      author={Feldman, I},
      author={Gohberg, I.},
      author={Krupnik, N.},
       title={{A method of explicit factorization of matrix functions and
  applications}},
        date={1994},
     journal={Integral Equations Operator Theory},
      volume={18},
      number={3},
       pages={277\ndash 302},
}

\bib{FGK95}{article}{
      author={Feldman, I.},
      author={Gohberg, I.},
      author={Krupnik, N.},
       title={{On explicit factorization and applications}},
        date={1995},
     journal={Integral Equations Operator Theory},
      volume={21},
       pages={430\ndash 459},
}

\bib{FGK00}{article}{
      author={Feldman, I.},
      author={Gohberg, I.},
      author={Krupnik, N.},
       title={{Convolution equations on finite intervals and factorization of
  matrix functions}},
        date={2000},
     journal={Integral Equations Operator Theory},
      volume={36},
      number={2},
       pages={201\ndash 211},
}

\bib{FGK04}{article}{
      author={Feldman, I.},
      author={Gohberg, I.},
      author={Krupnik, N.},
       title={{An explicit factorization algorithm}},
        date={2004},
     journal={Integral Equations Operator Theory},
      volume={49},
       pages={149\ndash 164},
}

\bib{FokSpen12}{article}{
      author={Fokas, A.~S.},
      author={Spence, E.~A.},
       title={Synthesis, as opposed to separation, of variables},
        date={2012},
        ISSN={0036-1445},
     journal={SIAM Rev.},
      volume={54},
      number={2},
       pages={291\ndash 324},
         url={https://doi-org.manchester.idm.oclc.org/10.1137/100809647},
      review={\MR{2916309}},
}

\bib{Fre03}{article}{
      author={Fredholm, I.},
       title={{Sur une classe d'\'{e}quations fonctionnelles} (french)},
        date={1903},
     journal={Acta Math.},
      volume={27},
       pages={365\ndash 390},
}

\bib{Gak50}{article}{
      author={Gakhov, F.~D.},
       title={{One case of the Riemann boundary value problem for a system of
  $n$ pairs of functions} (russian)},
        date={1950},
     journal={Izv. AN SSSR Ser. Mat.},
      volume={14},
       pages={549\ndash 568},
}

\bib{Gak52}{article}{
      author={Gakhov, F.~D.},
       title={{Riemann's boundary value problem for a system of n pairs of
  functions} (russian)},
        date={1952},
     journal={Izv. AN SSSR Ser. Mat.},
      volume={VII},
      number={4 (50)},
       pages={3\ndash 54},
}

\bib{Gak77}{book}{
      author={Gakhov, F.~D.},
       title={Boundary value problems},
     edition={3rd ed. (Russian)},
   publisher={Nauka},
     address={Moscow},
        date={1977},
}

\bib{Ga-Che}{book}{
      author={Gakhov, F.~D.},
      author={{\v{C}}erski{\u\i}, Ju.~I.},
       title={Equations of convolution type (in {R}ussian)},
   publisher={Nauka},
     address={Moscow},
        date={1978},
      review={\MR{527628 (80j:45001)}},
}

\bib{Spit}{incollection}{
      author={Gohberg, I.},
      author={Kaashoek, M.A.},
      author={Spitkovsky, I.M.},
       title={An overview of matrix factorization theory and operator
  applications},
    language={English},
        date={2003},
   booktitle={Factorization and integrable systems},
      editor={Gohberg, Israel},
      editor={Manojlovic, Nenad},
      editor={dos Santos, Ant},
      series={Operator Theory: Advances and Applications},
      volume={141},
   publisher={Birkh\"auser Basel},
       pages={1\ndash 102},
         url={http://dx.doi.org/10.1007/978-3-0348-8003-9_1},
}

\bib{GohKre58}{article}{
      author={Gohberg, I.},
      author={Krein, M.~G.},
       title={{Systems of integral equations on a half line with kernels
  depending on the difference of arguments}},
        date={1960},
     journal={Am. Math. Soc. Transl., Ser. 2},
      volume={24},
       pages={217\ndash 287},
}

\bib{Gohberg74}{book}{
      author={Gohberg, I.~C.},
      author={Fel\cprime~dman, I.~A.},
       title={Convolution equations and projection methods for their solution},
   publisher={American Mathematical Society, Providence, R.I.},
        date={1974},
        note={Translated from the Russian by F. M. Goldware, Translations of
  Mathematical Monographs, Vol. 41},
      review={\MR{0355675}},
}

\bib{GT}{article}{
      author={Gorbushin, N.},
      author={Truskinovsky, L.},
       title={Peristalsis by pulses of activity},
        date={2021Apr},
     journal={Phys. Rev. E},
      volume={103},
       pages={042411},
         url={https://link.aps.org/doi/10.1103/PhysRevE.103.042411},
}

\bib{GMT}{article}{
      author={Gorbushin, Nikolai},
      author={Eremeyev, Victor~A.},
      author={Mishuris, Gennady},
       title={On stress singularity near the tip of a crack with surface
  stresses},
        date={2020},
        ISSN={0020-7225},
     journal={Internat. J. Engrg. Sci.},
      volume={146},
       pages={103183, 17},
  url={https://doi-org.manchester.idm.oclc.org/10.1016/j.ijengsci.2019.103183},
      review={\MR{4023137}},
}

\bib{Gower_19}{article}{
      author={Gower, Artur~L.},
      author={Abrahams, David},
      author={Parnell, William~J.},
       title={A proof that multiple waves propagate in ensemble-averaged
  particulate materials},
        date={2019},
        ISSN={1364-5021},
     journal={Proc. A.},
      volume={475},
      number={2229},
       pages={20190344, 21},
         url={https://doi-org.manchester.idm.oclc.org/10.1098/rspa.2019.0344},
      review={\MR{4017266}},
}

\bib{finance}{article}{
      author={Green, Ross},
      author={Fusai, Gianluca},
      author={Abrahams, I~David},
       title={The {W}iener-{H}opf technique and discretely monitored
  path-dependent option pricing},
        date={2010},
        ISSN={0960-1627},
     journal={Math. Finance},
      volume={20},
      number={2},
       pages={259\ndash 288},
         url={http://dx.doi.org/10.1111/j.1467-9965.2010.00397.x},
      review={\MR{2650248 (2011h:91201)}},
}

\bib{Grudsky_08}{article}{
      author={Grudsky, Sergei},
      author={Shargorodsky, Eugene},
       title={Spectra of {T}oeplitz operators and compositions of {M}uckenhoupt
  weights with {B}laschke products},
        date={2008},
        ISSN={0378-620X},
     journal={Integral Equations Operator Theory},
      volume={61},
      number={1},
       pages={63\ndash 75},
  url={https://doi-org.manchester.idm.oclc.org/10.1007/s00020-008-1583-8},
      review={\MR{2414440}},
}

\bib{Mov_Cras_17}{article}{
      author={Haslinger, S.~G.},
      author={Movchan, N.~V.},
      author={Movchan, A.~B.},
      author={Jones, I.~S.},
      author={Craster, R.~V.},
       title={Controlling flexural waves in semi-infinite platonic crystals
  with resonator-type scatterers},
        date={2017},
        ISSN={0033-5614},
     journal={Quart. J. Mech. Appl. Math.},
      volume={70},
      number={3},
       pages={215\ndash 247},
         url={https://doi-org.manchester.idm.oclc.org/10.1093/qjmam/hbx005},
      review={\MR{3712202}},
}

\bib{Nigel_poles}{article}{
      author={Heaton, C.~J.},
      author={Peake, N.},
       title={Acoustic scattering in a duct with mean swirling flow},
        date={2005},
        ISSN={0022-1120},
     journal={J. Fluid Mech.},
      volume={540},
       pages={189\ndash 220},
         url={http://dx.doi.org/10.1017/S0022112005005719},
      review={\MR{2263126 (2007e:76231)}},
}

\bib{Hei50}{article}{
      author={Heins, A.~E.},
       title={{Systems of Wiener-Hopf equations}},
        date={1960},
     journal={Am. Math. Soc. Transl., Ser. 2},
      volume={24},
       pages={217\ndash 287},
}

\bib{Hil12}{book}{
      author={Hilbert, D.},
       title={Grundz\"uge der integralgleichungen, drittes abschnitt.
  (german)},
   publisher={Teubner},
     address={Leipzig-Berlin},
        date={1912},
}

\bib{Hop34}{book}{
      author={Hopf, E.},
       title={Mathematical problems of radiative equilibrium},
   publisher={Cambridge Univ. Press},
     address={Cambridge},
        date={1934},
}

\bib{Hur76}{article}{
      author={Hurd, R.~A.},
       title={{The Wiener--Hopf--Hilbert method for diffraction problems}},
        date={1976},
     journal={Can. J. Phys.},
      volume={54},
      number={7},
       pages={775\ndash 780},
}

\bib{Hur87}{article}{
      author={Hurd, R.~A.},
       title={{The explicit factorization of $2 \times 2$ Wiener-Hopf
  matrices}},
        date={1987},
     journal={Technische Hochschule Darmstadt, Preprint},
      volume={-},
      number={1040},
}

\bib{owl}{article}{
      author={Jaworski, Justin~W.},
      author={Peake, N.},
       title={Aerodynamic noise from a poroelastic edge with implications for
  the silent flight of owls},
        date={2013},
        ISSN={0022-1120},
     journal={J. Fluid Mech.},
      volume={723},
       pages={456\ndash 479},
  url={http://0-dx.doi.org.pugwash.lib.warwick.ac.uk/10.1017/jfm.2013.139},
      review={\MR{3045100}},
}

\bib{Jon84}{article}{
      author={Jones, D.~S.},
       title={{Commutative Wiener-Hopf factorization of a matrix}},
        date={1984},
     journal={R. Soc. Lond. Proc. Ser. A Math. Phys. Eng. Sci.},
      volume={393},
       pages={185\ndash 192},
}

\bib{Jones84}{article}{
      author={Jones, D.~S.},
       title={Factorization of a {W}iener-{H}opf matrix},
        date={1984},
        ISSN={0272-4960},
     journal={IMA J. Appl. Math.},
      volume={32},
      number={1-3},
       pages={211\ndash 220},
         url={http://dx.doi.org/10.1093/imamat/32.1-3.211},
      review={\MR{740458 (85j:45012)}},
}

\bib{Khr71a}{article}{
      author={Khrapkov, A.~A.},
       title={{Certain cases of the elastic equilibrium of an infinite wedge
  with a non-symmetric notch at the vertex, subjected to concentrated forces}},
        date={1971},
     journal={J. Appl. Math. Mech. (PMM)},
      volume={35},
       pages={625\ndash 637},
}

\bib{Khr71b}{article}{
      author={Khrapkov, A.~A.},
       title={{Closed form solutions of problems on the elastic equilibrium of
  an infinite wedge with nonsymmetric notch at the apex}},
        date={1971},
     journal={J. Appl. Math. Mech. (PMM)},
      volume={35},
       pages={1009\ndash 1016},
}

\bib{Khr01}{book}{
      author={Khrapkov, A.~A.},
       title={Wiener-hopf method in mixed elasticity problems},
   publisher={VNIIG Publishing House},
     address={Sankt Peterburg},
        date={2001},
}

\bib{my+Lorna}{article}{
      author={Kisil, A.},
      author={Ayton, L.~J.},
       title={Aerodynamic noise from rigid trailing edges with finite porous
  extensions},
        date={2018},
     journal={Journal of Fluid Mechanics},
      volume={836},
}

\bib{My1}{article}{
      author={Kisil, Anastasia~V.},
       title={A constructive method for an approximate solution to scalar
  {W}iener--{H}opf equations},
        date={2013},
        ISSN={1364-5021},
     journal={Proc. R. Soc. Lond. Ser. A Math. Phys. Eng. Sci.},
      volume={469},
      number={2154},
       pages={20120721, 17},
  url={http://0-dx.doi.org.pugwash.lib.warwick.ac.uk/10.1098/rspa.2012.0721},
      review={\MR{3035429}},
}

\bib{mythesis}{thesis}{
      author={Kisil, Anastasia~V.},
       title={Approximate {W}iener--{H}opf factorisation with stability
  analysis},
        type={Ph.D. Thesis},
        date={2015},
}

\bib{MyStrip}{article}{
      author={Kisil, Anastasia~V.},
       title={The relationship between a strip {W}iener-{H}opf problem and a
  line {R}iemann-{H}ilbert problem},
        date={2015},
        ISSN={0272-4960},
     journal={IMA J. Appl. Math.},
      volume={80},
      number={5},
       pages={1569\ndash 1581},
         url={https://doi-org.manchester.idm.oclc.org/10.1093/imamat/hxv007},
      review={\MR{3403991}},
}

\bib{MyD-K}{article}{
      author={Kisil, Anastasia~V.},
       title={Stability analysis of matrix {W}iener-{H}opf factorization of
  {D}aniele-{K}hrapkov class and reliable approximate factorization},
        date={2015},
        ISSN={1364-5021},
     journal={Proc. A.},
      volume={471},
      number={2178},
       pages={20150146, 15},
         url={https://doi-org.manchester.idm.oclc.org/10.1098/rspa.2015.0146},
      review={\MR{3367316}},
}

\bib{My_iter}{article}{
      author={Kisil, Anastasia~V.},
       title={An iterative {W}iener--{H}opf method for triangular matrix
  functions with exponential factors},
        date={2018},
     journal={SIAM Journal on Applied Mathematics},
      volume={78},
      number={1},
       pages={45\ndash 62},
      eprint={https://doi.org/10.1137/17M1136304},
         url={https://doi.org/10.1137/17M1136304},
}

\bib{Kiy12}{article}{
      author={Kiyasov, S.~N.},
       title={{Some cases of efficient factorization of second-order matrix
  functions}},
        date={2012},
     journal={Russian Math. (Iz. VUZ)},
      volume={56},
      number={6},
       pages={36\ndash 43},
}

\bib{Kiy13}{article}{
      author={Kiyasov, S.~N.},
       title={{Certain classes of problems on linear conjugation for a
  two-dimensional vector admitting explicit solution}},
        date={2013},
     journal={Russian Math. (Iz. VUZ)},
      volume={57},
      number={1},
       pages={1\ndash 16},
}

\bib{Kiy15}{article}{
      author={Kiyasov, S.~N.},
       title={{Some classes of linear conjugation problems for a
  three-dimensional vector that are solvable in closed form}},
        date={2015},
     journal={Siberian Math. J.},
      volume={56},
      number={2},
       pages={313\ndash 329},
}

\bib{Kobayashi91}{article}{
      author={Kobayashi, Kazuya},
       title={Plane wave diffraction by a strip: exact and asymptotic
  solutions},
        date={1991},
        ISSN={0031-9015},
     journal={J. Phys. Soc. Japan},
      volume={60},
      number={6},
       pages={1891\ndash 1905},
         url={https://doi-org.manchester.idm.oclc.org/10.1143/JPSJ.60.1891},
      review={\MR{1121835}},
}

\bib{Koiter}{article}{
      author={Koiter, W.~T.},
       title={Approximate solution of {W}iener-{H}opf type integral equations
  with applications. {I}. {G}eneral theory},
        date={1954},
        ISSN={0023-3366},
     journal={Nederl. Akad. Wetensch. Proc. Ser. B.},
      volume={57},
       pages={558\ndash 564},
      review={\MR{0073854 (17,498e)}},
}

\bib{Kopev08}{article}{
      author={{Kop'ev}, V.~F.},
      author={{Faranosov}, G.~A.},
       title={{Control over the instability wave in terms of the
  two-dimensional model of a nozzle edge}},
        date={2008-05},
     journal={Acoustical Physics},
      volume={54},
      number={3},
       pages={319\ndash 326},
}

\bib{Kranzer_68}{article}{
      author={Kranzer, Herbert~C.},
       title={Asymptotic factorization in nondissipative {W}iener-{H}opf
  problems},
        date={1967},
     journal={J. Math. Mech.},
      volume={17},
       pages={577\ndash 600},
      review={\MR{0220020 (36 \#3087)}},
}

\bib{sto-W-H-eng}{article}{
      author={Kuznetsov, A.},
      author={Kyprianou, A.~E.},
      author={Pardo, J.~C.},
      author={van Schaik, K.},
       title={A {W}iener-{H}opf {M}onte {C}arlo simulation technique for
  {L}\'evy processes},
        date={2011},
        ISSN={1050-5164},
     journal={Ann. Appl. Probab.},
      volume={21},
      number={6},
       pages={2171\ndash 2190},
         url={http://dx.doi.org/10.1214/10-AAP746},
      review={\MR{2895413 (2012m:60109)}},
}

\bib{history}{article}{
      author={Lawrie, Jane~B.},
      author={Abrahams, I.~David},
       title={A brief historical perspective of the {W}iener-{H}opf technique},
        date={2007},
        ISSN={0022-0833},
     journal={J. Engrg. Math.},
      volume={59},
      number={4},
       pages={351\ndash 358},
         url={http://dx.doi.org/10.1007/s10665-007-9195-x},
      review={\MR{2373797 (2009a:45002)}},
}

\bib{Lighthill}{article}{
      author={Lighthill, M.~J.},
       title={On sound generated aerodynamically {I}. general theory},
        date={1952},
        ISSN={0080-4630},
     journal={Proc. R. Soc. Lond. Ser. A Math. Phys. Eng. Sci.},
      volume={211},
      number={1107},
       pages={564\ndash 587},
}

\bib{Handbook_water}{book}{
      author={Linton, C.~M.},
      author={McIver, P.},
       title={Handbook of mathematical techniques for wave/structure
  interactions},
   publisher={Chapman \& Hall/CRC, Boca Raton, FL},
        date={2001},
        ISBN={1-58488-132-1},
         url={https://doi.org/10.1201/9781420036060},
      review={\MR{1886390}},
}

\bib{Lit_shift}{book}{
      author={Litvinchuk, Georgii~S.},
       title={Solvability theory of boundary value problems and singular
  integral equations with shift},
      series={Mathematics and its Applications},
   publisher={Kluwer Academic Publishers, Dordrecht},
        date={2000},
      volume={523},
        ISBN={0-7923-6549-6},
  url={https://doi-org.manchester.idm.oclc.org/10.1007/978-94-011-4363-9},
      review={\MR{1885495}},
}

\bib{Spit_book}{book}{
      author={Litvinchuk, Georgii~S.},
      author={Spitkovskii, Ilia~M.},
       title={Factorization of measurable matrix functions},
      series={Operator Theory: Advances and Applications},
   publisher={Birkh\"auser Verlag, Basel},
        date={1987},
      volume={25},
        ISBN={3-7643-1883-X},
         url={http://dx.doi.org/10.1007/978-3-0348-6266-0},
        note={Translated from the Russian by Bernd Luderer, With a foreword by
  Bernd Silbermann},
      review={\MR{1015716 (90g:47030)}},
}

\bib{Mish_crack_iter}{article}{
      author={Livasov, P.},
      author={Mishuris, G.},
       title={Numerical factorization of a matrix-function with exponential
  factors in an anti-plane problem for a crack with process zone},
        date={2019},
     journal={Philosophical Transactions of the Royal Society A: Mathematical,
  Physical and Engineering Sciences},
      volume={377},
      number={2156},
       pages={20190109},
  eprint={https://royalsocietypublishing.org/doi/pdf/10.1098/rsta.2019.0109},
  url={https://royalsocietypublishing.org/doi/abs/10.1098/rsta.2019.0109},
}

\bib{nano_Wiener_19}{article}{
      author={Margetis, Dionisios},
       title={Edge plasmon-polaritons on isotropic semi-infinite conducting
  sheets},
        date={2020},
        ISSN={0022-2488},
     journal={J. Math. Phys.},
      volume={61},
      number={6},
       pages={062901, 25},
         url={https://doi-org.manchester.idm.oclc.org/10.1063/1.5128895},
      review={\MR{4112646}},
}

\bib{Scatterers_15}{article}{
      author={Martin, P.~A.},
      author={Abrahams, I.~David},
      author={Parnell, William~J.},
       title={One-dimensional reflection by a semi-infinite periodic row of
  scatterers},
        date={2015},
        ISSN={0165-2125},
     journal={Wave Motion},
      volume={58},
       pages={1\ndash 12},
  url={https://doi-org.manchester.idm.oclc.org/10.1016/j.wavemoti.2015.06.005},
      review={\MR{3378920}},
}

\bib{Maurya19}{article}{
      author={Maurya, Gaurav},
      author={Sharma, Basant~Lal},
       title={Scattering by two staggered semi-infinite cracks on square
  lattice: an application of asymptotic {W}iener-{H}opf factorization},
        date={2019},
        ISSN={0044-2275},
     journal={Z. Angew. Math. Phys.},
      volume={70},
      number={5},
       pages={Paper No. 133, 21},
  url={https://doi-org.manchester.idm.oclc.org/10.1007/s00033-019-1183-2},
      review={\MR{3999332}},
}

\bib{Mikhlin86}{book}{
      author={Mikhlin, Solomon~G.},
      author={Pr\"{o}ssdorf, Siegfried},
       title={Singular integral operators},
   publisher={Springer-Verlag, Berlin},
        date={1986},
        ISBN={3-540-15967-3},
  url={https://doi-org.manchester.idm.oclc.org/10.1007/978-3-642-61631-0},
        note={Translated from the German by Albrecht B\"{o}ttcher and Reinhard
  Lehmann},
      review={\MR{881386}},
}

\bib{Mishuris-Rogosin}{article}{
      author={Mishuris, G.},
      author={Rogosin, S.},
       title={An asymptotic method of factorization of a class of matrix
  functions},
        date={2014},
     journal={Proc. R. Soc. Lond. Ser. A Math. Phys. Eng. Sci.},
      volume={470},
      number={2166},
}

\bib{assym_2}{article}{
      author={Mishuris, G.},
      author={Rogosin, S.},
       title={Factorization of a class of matrix-functions with stable partial
  indices},
        date={2016},
        ISSN={0170-4214},
     journal={Math. Methods Appl. Sci.},
      volume={39},
      number={13},
       pages={3791\ndash 3807},
         url={http://dx.doi.org/10.1002/mma.3825},
      review={\MR{3529384}},
}

\bib{MR_2018}{article}{
      author={Mishuris, G.},
      author={Rogosin, S.},
       title={Regular approximate factorization of a class of matrix-function
  with an unstable set of partial indices},
        date={2018},
        ISSN={1364-5021},
     journal={Proc. A.},
      volume={474},
      number={2209},
       pages={20170279, 18},
         url={https://doi-org.manchester.idm.oclc.org/10.1098/rspa.2017.0279},
      review={\MR{3762905}},
}

\bib{Mishris96}{article}{
      author={Mishuris, G.~S.},
       title={Boundary value problems for {P}oisson's equation in a
  multi-wedge--multi-layered region. {II}. {G}eneral type of interfacial
  conditions},
        date={1996},
        ISSN={0373-2029},
     journal={Arch. Mech. (Arch. Mech. Stos.)},
      volume={48},
      number={4},
       pages={711\ndash 745},
      review={\MR{1416158}},
}

\bib{Mishuris97}{article}{
      author={Mishuris, G.~S.},
       title={{$2$}-{D} boundary value problems of thermoelasticity in a
  multi-wedge--multi-layered region. {II}. {S}ystems of integral equations},
        date={1997},
        ISSN={0373-2029},
     journal={Arch. Mech. (Arch. Mech. Stos.)},
      volume={49},
      number={6},
       pages={1135\ndash 1165},
      review={\MR{1613721}},
}

\bib{Mishuris01}{article}{
      author={Mishuris, G.~S.},
      author={Kuhn, G.},
       title={Asymptotic behaviour of the elastic solution near the tip of a
  crack situated at a nonideal interface},
        date={2001},
        ISSN={0044-2267},
     journal={ZAMM Z. Angew. Math. Mech.},
      volume={81},
      number={12},
       pages={811\ndash 826},
  url={https://doi-org.manchester.idm.oclc.org/10.1002/1521-4001(200112)81:12<811::AID-ZAMM811>3.0.CO;2-I},
      review={\MR{1872768}},
}

\bib{Mishuris99}{article}{
      author={Mishuris, Gennady~S.},
       title={Stress singularity at a crack tip for various intermediate zones
  in bimaterial structures (mode {III})},
        date={1999},
        ISSN={0020-7683},
     journal={Internat. J. Solids Structures},
      volume={36},
      number={7},
       pages={999\ndash 1015},
  url={https://doi-org.manchester.idm.oclc.org/10.1016/S0020-7683(97)00342-9},
      review={\MR{1666101}},
}

\bib{Mishuris09}{article}{
      author={Mishuris, Gennady~S.},
      author={Movchan, Alexander~B.},
      author={Slepyan, Leonid~I.},
       title={Dynamics of a bridged crack in a discrete lattice},
        date={2008},
        ISSN={0033-5614},
     journal={Quart. J. Mech. Appl. Math.},
      volume={61},
      number={2},
       pages={151\ndash 160},
         url={http://dx.doi.org/10.1093/qjmam/hbm030},
      review={\MR{2414430 (2009d:74085)}},
}

\bib{MishurisOlesiak95}{article}{
      author={Mishuris, Gennady~S.},
      author={Olesiak, Zbigniew~S.},
       title={On boundary value problems in fracture of elastic composites},
        date={1995},
        ISSN={0956-7925},
     journal={European J. Appl. Math.},
      volume={6},
      number={6},
       pages={591\ndash 610},
  url={https://doi-org.manchester.idm.oclc.org/10.1017/S0956792500002084},
      review={\MR{1364719}},
}

\bib{Moi89}{article}{
      author={Moiseev, N.~G.},
       title={{On factorization of matrix-valued functions of special form}},
        date={1989},
     journal={Sov. Math. Dokl.},
      volume={39},
      number={2},
       pages={264\ndash 267},
}

\bib{Moi93}{article}{
      author={Moiseev, N.~G.},
       title={{Partial factorization indices of a special class of
  matrix-valued functions}},
        date={1994},
     journal={Dokl. Akad. Nauk},
      volume={48},
      number={3},
       pages={640\ndash 644},
}

\bib{Muskh_book}{book}{
      author={Muskhelishvili, N.~I.},
       title={Singular integral equations. {B}oundary problems of function
  theory and their application to mathematical physics},
   publisher={P. Noordhoff N. V., Groningen},
        date={1953},
        note={Translation by J. R. M. Radok},
      review={\MR{0058845}},
}

\bib{Mus68}{book}{
      author={Muskhelishvili, N.~I.},
       title={Singular integral equations},
     edition={3rd ed. (Russian)},
   publisher={Nauka},
     address={Moscow},
        date={1968},
}

\bib{AAA_Trefethen}{article}{
      author={Nakatsukasa, Yuji},
      author={Se~te, Olivier},
      author={Trefethen, Lloyd~N.},
       title={The aaa algorithm for rational approximation},
    language={eng},
        date={2018},
        ISSN={1064-8275},
     journal={SIAM Journal on Scientific Computing},
      volume={40},
      number={3},
       pages={A1494\ndash A1522},
}

\bib{Nethercote2020a}{article}{
      author={Nethercote, M.},
      author={Assier, R.C.},
      author={Abrahams, I.D.},
       title={Analytical methods for perfect wedge diffraction: A review},
        date={2020-03},
     journal={Wave Motion},
      volume={93},
       pages={102479},
         url={https://doi.org/10.1016/j.wavemoti.2019.102479},
}

\bib{bookWH}{book}{
      author={Noble, B.},
       title={Methods based on the {W}iener-{H}opf technique for the solution
  of partial differential equations},
      series={International Series of Monographs on Pure and Applied
  Mathematics. Vol. 7},
   publisher={Pergamon Press},
     address={New York},
        date={1958},
      review={\MR{0102719 (21 \#1505)}},
}

\bib{Nowak09}{article}{
      author={Nowak, Micha\l~A.},
       title={Approximation methods for a class of discrete {W}iener-{H}opf
  equations},
        date={2009},
        ISSN={1232-9274},
     journal={Opuscula Math.},
      volume={29},
      number={3},
       pages={271\ndash 288},
  url={https://doi-org.manchester.idm.oclc.org/10.7494/OpMath.2009.29.3.271},
      review={\MR{2551900}},
}

\bib{olver_computing_2011}{article}{
      author={Olver, Sheehan},
       title={Computing the {Hilbert} transform and its inverse},
    language={en},
        date={2011-02},
        ISSN={0025-5718, 1088-6842},
     journal={Math. Comp.},
      volume={80},
       pages={1745\ndash 1767},
}

\bib{olver_change}{article}{
      author={Olver, Sheehan},
       title={Change of variable formulas for regularizing slowly decaying and
  oscillatory cauchy and hilbert transforms},
        date={2014},
     journal={Analysis and Applications},
      volume={12},
      number={04},
       pages={369\ndash 384},
      eprint={https://doi.org/10.1142/S0219530514500353},
         url={https://doi.org/10.1142/S0219530514500353},
}

\bib{Piccolroaz07}{article}{
      author={Piccolroaz, A.},
      author={Mishuris, G.},
      author={Movchan, A.~B.},
       title={Evaluation of the {L}azarus-{L}eblond constants in the asymptotic
  model of the interfacial wavy crack},
        date={2007},
        ISSN={0022-5096},
     journal={J. Mech. Phys. Solids},
      volume={55},
      number={8},
       pages={1575\ndash 1600},
  url={https://doi-org.manchester.idm.oclc.org/10.1016/j.jmps.2007.02.001},
      review={\MR{2349279}},
}

\bib{Pricc09}{article}{
      author={Piccolroaz, A.},
      author={Mishuris, G.},
      author={Movchan, A.~B.},
       title={Symmetric and skew-symmetric weight functions in 2{D}
  perturbation models for semi-infinite interfacial cracks},
        date={2009},
        ISSN={0022-5096},
     journal={J. Mech. Phys. Solids},
      volume={57},
      number={9},
       pages={1657\ndash 1682},
  url={https://doi-org.manchester.idm.oclc.org/10.1016/j.jmps.2009.05.003},
      review={\MR{2555190}},
}

\bib{Ple08}{article}{
      author={Plemelj, J.},
       title={{Riemannsche Funktionenscharen mit gegebener Monodromiegruppe}},
        date={1908},
     journal={Monat. Math. Phys. (German)},
      volume={19},
       pages={211\ndash 245},
}

\bib{iter_n}{article}{
      author={Priddin, Matthew~J.},
      author={Kisil, Anastasia~V.},
      author={Ayton, Lorna~J.},
       title={Applying an iterative method numerically to solve n by n matrix
  wiener--hopf equations with exponential factors},
        date={2020},
     journal={Philosophical Transactions of the Royal Society A: Mathematical,
  Physical and Engineering Sciences},
      volume={378},
      number={2162},
       pages={20190241},
  eprint={https://royalsocietypublishing.org/doi/pdf/10.1098/rsta.2019.0241},
  url={https://royalsocietypublishing.org/doi/abs/10.1098/rsta.2019.0241},
}

\bib{PriRog18}{article}{
      author={Primachuk, L.},
      author={Rogosin, S.},
       title={{Factorization of triangular matrix-functions of an arbitrary
  order}},
        date={2018},
     journal={Lobachevskii J. Math.},
      volume={39},
      number={1},
       pages={129\ndash 137},
}

\bib{Primachuk_20}{article}{
      author={Primachuk, L.},
      author={Rogosin, S.},
      author={Dubatovskaya, M.},
       title={On {$\Bbb R$}-linear conjugation problem on the unit circle},
        date={2020},
        ISSN={2077-9879},
     journal={Eurasian Math. J.},
      volume={11},
      number={3},
       pages={79\ndash 88},
  url={https://doi-org.manchester.idm.oclc.org/10.32523/2077-9879-2020-11-3-79-88},
      review={\MR{4162274}},
}

\bib{Prossdorf91}{book}{
      author={Pr\"{o}ssdorf, Siegfried},
      author={Silbermann, Bernd},
       title={Numerical analysis for integral and related operator equations},
      series={Operator Theory: Advances and Applications},
   publisher={Birkh\"{a}user Verlag, Basel},
        date={1991},
      volume={52},
        ISBN={3-7643-2620-4},
      review={\MR{1193030}},
}

\bib{Application_DK}{article}{
      author={Rawlins, A.~D.},
       title={Two waveguide trifurcation problems},
        date={1997},
        ISSN={0305-0041},
     journal={Math. Proc. Cambridge Philos. Soc.},
      volume={121},
      number={3},
       pages={555\ndash 573},
         url={http://dx.doi.org/10.1017/S0305004196001296},
      review={\MR{1434661 (97i:76117)}},
}

\bib{Raw}{article}{
      author={Rawlins, A.~D.},
      author={Williams, W.~E.},
       title={Matrix {W}iener-{H}opf factorisation},
        date={1981},
        ISSN={0033-5614},
     journal={Quart. J. Mech. Appl. Math.},
      volume={34},
      number={1},
       pages={1\ndash 8},
         url={http://dx.doi.org/10.1093/qjmam/34.1.1},
      review={\MR{610745 (82e:15020)}},
}

\bib{Rogosin_21}{article}{
      author={Rogosin, S},
      author={Khvoshchinskaya, L},
       title={Construction of ${2\times 2}$ fuchsian system with five singular
  points},
    language={eng},
        date={2021},
        ISSN={1995-0802},
     journal={Lobachevskii journal of mathematics},
      volume={42},
      number={4},
       pages={830\ndash 849},
}

\bib{RogMis16}{article}{
      author={Rogosin, S.},
      author={Mishuris, G.},
       title={{Constructive methods for factorization of matrix-functions}},
        date={2016},
     journal={IMA Journal of Applied Mathematics},
      volume={81},
      number={2},
       pages={365\ndash 391},
}

\bib{RogDub17}{inproceedings}{
      author={Rogosin, S.~V.},
      author={Dubatovskaya, M.~V.},
       title={On solution to the {R}iemann vector-matrix boundary value problem
  with funcionally-commutative coefficient},
        date={2017},
   booktitle={Proc. intern. conf., october 10-12, 2017},
      editor={et~al., A.A.~Bolshakov},
      series={MMTT},
      volume={10},
   publisher={St.-Petersburg Polytechn. University},
     address={St.-Petersburg},
       pages={11\ndash 15},
}

\bib{Sakhnovich_15}{book}{
      author={Sakhnovich, L.~A.},
       title={Integral equations with difference kernels on finite intervals},
    language={eng},
     edition={Second edition, revised and extended.},
      series={Operator theory: advances and applications, volume 84},
   publisher={Birkhauser},
     address={Cham},
        date={2015},
        ISBN={9783319164892},
}

\bib{Sautbekov_17}{article}{
      author={Sautbekova, Merey},
      author={Sautbekov, Seil},
       title={Solution of the {N}eumann problem of diffraction by a strip using
  the {W}iener-{H}opf method: short-wave asymptotic solutions},
        date={2017},
        ISSN={0018-926X},
     journal={IEEE Trans. Antennas and Propagation},
      volume={65},
      number={9},
       pages={4797\ndash 4802},
  url={https://doi-org.manchester.idm.oclc.org/10.1109/TAP.2017.2730923},
      review={\MR{3697000}},
}

\bib{Schwarzschild1901}{article}{
      author={Schwarzschild, K.},
       title={Die beugung und polarisation des lichts durch einen spalt. {I}},
        date={1901},
        ISSN={1432-1807},
     journal={Mathematische Annalen},
      volume={55},
      number={2},
       pages={177\ndash 247},
         url={http://dx.doi.org/10.1007/BF01444971},
}

\bib{Shanin_13}{article}{
      author={Shanin, A.~V.},
       title={Solution of {R}iemann-{H}ilbert problem related to
  {W}iener-{H}opf matrix factorization problem using ordinary differential
  equations in the commutative case},
        date={2013},
        ISSN={0033-5614},
     journal={Quart. J. Mech. Appl. Math.},
      volume={66},
      number={4},
       pages={533\ndash 555},
         url={http://dx.doi.org/10.1093/qjmam/hbt017},
      review={\MR{3138009}},
}

\bib{Shanin20}{article}{
      author={Shanin, A.~V.},
      author={Korolkov, A.~I.},
       title={Sommerfeld-type integrals for discrete diffraction problems},
        date={2020},
        ISSN={0165-2125},
     journal={Wave Motion},
      volume={97},
       pages={102606, 24},
  url={https://doi-org.manchester.idm.oclc.org/10.1016/j.wavemoti.2020.102606},
      review={\MR{4107117}},
}

\bib{Shanin_S_3}{article}{
      author={Shanin, Andrey~V.},
       title={Three theorems concerning diffraction by a strip or a slit},
        date={2001},
     journal={The Quarterly Journal of Mechanics and Applied Mathematics},
      volume={54},
      number={1},
       pages={107},
  eprint={/oup/backfile/content_public/journal/qjmam/54/1/10.1093/qjmam/54.1.107/2/540107.pdf},
         url={+ http://dx.doi.org/10.1093/qjmam/54.1.107},
}

\bib{comm}{unpublished}{
      author={Shanin, A.V.},
      author={Doubravsky, E.M.},
       title={Criteria for commutative factorization of a class of algebraic
  matrices},
        note={http://arxiv.org/abs/1211.4424},
}

\bib{SHARMA17_Somm}{article}{
      author={Sharma, Basant~Lal},
       title={Continuum limit of discrete {S}ommerfeld problems on square
  lattice},
        date={2017},
        ISSN={0256-2499},
     journal={S\={a}dhan\={a}},
      volume={42},
      number={5},
       pages={713\ndash 728},
      review={\MR{3659067}},
}

\bib{Sharma20}{article}{
      author={Sharma, Basant~Lal},
      author={Mishuris, Gennady},
       title={Scattering on a square lattice from a crack with a damage zone},
        date={2020},
        ISSN={1364-5021},
     journal={Proc. A.},
      volume={476},
      number={2235},
       pages={20190686, 17},
         url={https://doi-org.manchester.idm.oclc.org/10.1137/15m1010646},
      review={\MR{4089002}},
}

\bib{Slepyan15}{article}{
      author={Slepyan, Leonid},
      author={Ayzenberg-Stepanenko, Mark},
      author={Mishuris, Gennady},
       title={Forerunning mode transition in a continuous waveguide},
        date={2015},
        ISSN={0022-5096},
     journal={J. Mech. Phys. Solids},
      volume={78},
       pages={32\ndash 45},
  url={https://doi-org.manchester.idm.oclc.org/10.1016/j.jmps.2015.01.015},
      review={\MR{3349449}},
}

\bib{Slepyan02}{book}{
      author={Slepyan, Leonid~I.},
       title={Models and phenomena in fracture mechanics},
      series={Foundations of Engineering Mechanics},
   publisher={Springer-Verlag, Berlin},
        date={2002},
        ISBN={3-540-43767-3},
  url={https://doi-org.manchester.idm.oclc.org/10.1007/978-3-540-48010-5},
      review={\MR{1986072}},
}

\bib{Smith_20}{article}{
      author={Smith, M. J.~A.},
      author={Peter, M.~A.},
      author={Abrahams, I.~D.},
      author={Meylan, M.~H.},
       title={On the {W}iener-{H}opf solution of water-wave interaction with a
  submerged elastic or poroelastic plate},
        date={2020},
        ISSN={1364-5021},
     journal={Proc. A.},
      volume={476},
      number={2242},
       pages={20200360, 21},
      review={\MR{4176311}},
}

\bib{Sneddon95}{book}{
      author={Sneddon, Ian~N.},
       title={Fourier transforms},
   publisher={Dover Publications, Inc., New York},
        date={1995},
        ISBN={0-486-68522-5},
        note={Reprint of the 1951 original},
      review={\MR{1348163}},
}

\bib{Speck_19}{article}{
      author={Speck, Frank-Olme},
       title={From {S}ommerfeld diffraction problems to operator
  factorisation},
        date={201912},
     journal={Constructive Mathematical Analysis},
      volume={2},
       pages={183\ndash 216},
}

\bib{Stahl97}{article}{
      author={Stahl, Herbert},
       title={The convergence of {P}ad\'{e} approximants to functions with
  branch points},
        date={1997},
        ISSN={0021-9045},
     journal={J. Approx. Theory},
      volume={91},
      number={2},
       pages={139\ndash 204},
         url={https://doi-org.manchester.idm.oclc.org/10.1006/jath.1997.3141},
      review={\MR{1484040}},
}

\bib{ice_free_surface}{article}{
      author={{Tkacheva}, L.~A.},
       title={{Wave Pattern Due to a Load Moving on the Free Surface of a Fluid
  along the Edge of an Ice Sheet}},
        date={2019-05},
     journal={Journal of Applied Mechanics and Technical Physics},
      volume={60},
       pages={462\ndash 472},
}

\bib{trogdon_riemannhilbert_2016}{book}{
      author={Trogdon, Tom},
      author={Olver, Sheehan},
       title={Riemann-{Hilbert} {Problems}, {Their} {Numerical} {Solution}, and
  the {Computation} of {Nonlinear} {Special} {Functions}},
    language={en},
   publisher={SIAM},
        date={2016},
}

\bib{Thompson_14}{article}{
      author={Tymis, N.},
      author={Thompson, I.},
       title={Scattering by a semi-infinite lattice and the excitation of
  {B}loch waves},
        date={2014},
        ISSN={0033-5614},
     journal={Quart. J. Mech. Appl. Math.},
      volume={67},
      number={3},
       pages={469\ndash 503},
         url={https://doi-org.manchester.idm.oclc.org/10.1093/qjmam/hbu014},
      review={\MR{3249489}},
}

\bib{ufimtsev2003theory}{book}{
      author={Ufimtsev, P.Y.},
      author={Terzuoli, A.J.},
      author={Moore, R.D.},
       title={Theory of edge diffraction in electromagnetics},
      series={Contemporary research on emerging sciences and technology},
   publisher={Tech Science Press},
        date={2003},
        ISBN={9780965700177},
         url={https://books.google.co.uk/books?id=lTVDAAAACAAJ},
}

\bib{Nigel_cyl}{article}{
      author={Veitch, B.},
      author={Peake, N.},
       title={Acoustic propagation and scattering in the exhaust flow from
  coaxial cylinders},
        date={2008},
        ISSN={0022-1120},
     journal={J. Fluid Mech.},
      volume={613},
       pages={275\ndash 307},
         url={http://dx.doi.org/10.1017/S0022112008003169},
      review={\MR{2462429 (2009j:76203)}},
}

\bib{VeiAbr07}{article}{
      author={Veitch, B.~H.},
      author={Abrahams, I.~D.},
       title={{On the commutative factorization of $n\times n$ matrix
  Wiener-Hopf kernels with distinct eigenvalues}},
        date={2007},
     journal={R. Soc. Lond. Proc. Ser. A Math. Phys. Eng. Sci.},
      volume={463},
       pages={613\ndash 639},
}

\bib{Vekua}{book}{
      author={Vekua, N.~P.},
       title={Systems of singular integral equations},
   publisher={P. Noordhoff, Ltd., Groningen},
        date={1967},
        note={Translated from the Russian by A. G. Gibbs and G. M. Simmons.
  Edited by J. H. Ferziger},
      review={\MR{0211220 (35 \#2102)}},
}

\bib{AAM_2021}{unpublished}{
      author={V.M.~Adukov, G.~Mishuris, N.V.~Adukova},
       title={An explicit wiener–hopf factorization of matrix
  polynomials:exact versus numerical solutions of the problem},
        note={Submitted Proc. R. Soc. A.},
}

\bib{Mag_edge_88}{article}{
      author={{Volkov}, V.~A.},
      author={{Mikhailov}, S.~A.},
       title={{Edge magnetoplasmons - Low-frequency weakly damped excitations
  in homogeneous two-dimensional electron systems}},
        date={1988-08},
     journal={Zhurnal Eksperimentalnoi i Teoreticheskoi Fiziki},
      volume={94},
       pages={217\ndash 241},
}

\bib{Vol09}{article}{
      author={Volterra, V.},
       title={{Sulle equa\-zioni integro-dif\-feren\-ziali}},
        date={1909},
     journal={Atti della Reale Accademia dei Lincei, Rendiconti Classe Sci.
  Fis. Mat e Nat. (Ser. V) (Italian)},
      volume={18},
       pages={203\ndash 211},
}

\bib{vorovich1979}{book}{
      author={Vorovich, I.I.},
      author={Babeshko, V.A.},
       title={Dynamical mixed boundary problems in elastisity for non-classical
  domains (in {R}ussian)},
   publisher={Nauka},
        date={1979},
}

\bib{WieHop31}{article}{
      author={Wiener, N.},
      author={Hopf, E.},
       title={{\"Uber eine Klasse singularer Integralgleichungen}},
        date={1931},
     journal={Sitzg.-ber. Preuss. Akad. Wiss. ; Phys.-Math. (German)},
      volume={-},
      number={30--32},
       pages={696\ndash 706},
}

\bib{Water_equi}{article}{
      author={Williams, T.~D.},
      author={Meylan, Michael~H.},
       title={The wiener--hopf and residue calculus solutions for a submerged
  semi-infinite elastic plate},
        date={2012Aug},
        ISSN={1573-2703},
     journal={Journal of Engineering Mathematics},
      volume={75},
      number={1},
       pages={81\ndash 106},
         url={https://doi.org/10.1007/s10665-011-9518-9},
}

\bib{Willis95}{article}{
      author={Willis, J.~R.},
      author={Movchan, A.~B.},
       title={Dynamic weight functions for a moving crack. {I}. {M}ode {I}
  loading},
        date={1995},
        ISSN={0022-5096},
     journal={J. Mech. Phys. Solids},
      volume={43},
      number={3},
       pages={319\ndash 341},
  url={https://doi-org.manchester.idm.oclc.org/10.1016/0022-5096(94)00075-G},
      review={\MR{1321106}},
}

\bib{old_iter}{article}{
      author={Wu, {Tai Te}},
       title={Iterative solution of {W}iener-{H}opf integral equation},
        date={1963},
     journal={Q Appl. Math.},
      volume={20},
       pages={341\ndash 352},
}

\bib{Pade_numeric_test}{article}{
      author={Yamada, Hiroaki},
      author={Ikeda, Kensuke},
       title={A numerical test of pade approximation for some functions with
  singularity},
    language={eng},
        date={2014},
     journal={arXiv.org},
         url={http://search.proquest.com/docview/2084158312/},
}

\bib{DK_plas}{article}{
      author={Yener, S.},
      author={Serbest, A.H.},
       title={Diffraction of plane waves by an impedance half-plane in cold
  plasma},
        date={2002},
     journal={Journal of Electromagnetic Waves and Applications},
      volume={16},
      number={7},
       pages={995\ndash 1005},
      eprint={https://doi.org/10.1163/156939302X00318},
         url={https://doi.org/10.1163/156939302X00318},
}

\bib{Zve71}{article}{
      author={Zverovich, E.~I.},
       title={{Boundary value problems in the theory of analytic functions in
  H\"older classes on Riemann surfaces}},
        date={1971},
     journal={Uspekhi Mat. Nauk (Russian)},
      volume={26},
      number={1 (157)},
       pages={113\ndash 179},
}

\end{biblist}
\end{bibdiv}

\end{document}